\documentclass[12pt,a4paper]{article}
\usepackage{amsthm}
\usepackage{fleqn}
\usepackage{array}
\usepackage[latin1]{inputenc}
\usepackage{amssymb}
\usepackage{amsfonts}

\newcommand{\bR}{{\mathbb{R}}}

\newcommand{\bZ}{{\mathbb{Z}}}

\newcommand{\n}{\nu}

\usepackage[T1]{fontenc}           
\usepackage{textcomp}              
\usepackage{epsfig}        

\theoremstyle{definition}
\newtheorem{definition}{Definition}[section]

\usepackage{showkeys}

\begin{document}
\begin{center}
    \textbf{\Large 
Density estimates for $k$-impassable lattices of balls and general convex
bodies in ${\Bbb R}^n$}\\[2ex]
\emph{E. Makai, Jr., \footnote{The author was partially supported by the
Hungarian National
Foundation for Scientific Research, Grant Nos. T046846, T043520, and K68398.}
and H. Martini\footnote{The author was supported by the
``Discrete and Convex'' Project (MTKD-CT-2005-014333)
carried out by the A. Rényi Institute of Mathematics, Hungarian Academy of
Sciences, in the framework of the
European Community's ``Structuring the European Research Area'' Program.}}
\end{center}

\begin{abstract}
G. Fejes Tóth posed the following problem: Determine the infimum
of the densities of the lattices of closed balls in $\bR^n$ such
that each affine $k$-subspace $(0 \le k \le n-1)$ of $\bR^n$ intersects some 
ball of
the lattice. 
We give a lower estimate
for any $n,k$ like above. If, in the problem posed by G. Fejes Tóth, we
replace the ball $B^n$ by a (centrally symmetric) convex body $K
\subset \bR^n$, we may ask for the infimum of all above infima of
densities of lattices of translates of $K$ with the above
property, when $K$ ranges over all (centrally symmetric) convex
bodies in $\bR^n$. For these quantities we give lower estimates as well,
which are sharp, or almost sharp, for certain classes of convex bodies $K$.
For $k=n-1$ we
give an upper estimate for the supremum of all above infima of
densities, $K$ also ranging as above (i.e., a ``minimax'' problem). For $n=2$ 
our
estimate
is rather close to the conjecturable maximum.
We point out the connection of the above questions to the following
problem: Find the largest radius of a cylinder, with base an
$(n-1)$-ball, that can
be fitted into any lattice packing of balls (actually, here balls can be
replaced by some convex bodies $K \subset \bR^n$, the axis of the cylinder
may be $k$-dimensional and its basis has to be
chosen suitably). Among others we complete the proof of a theorem
of I. Hortob\'agyi from 1971. Our proofs for the lower estimates of
densities for balls, and for the cylinder problem, follow quite
closely a paper of J. Horv\'ath from 1970. This paper is also an
addendum to a paper of the first named author 
from 1978 in the sense that to some
arguments given there not in a detailed manner, we give here for all of these
complete proofs.
\end{abstract}

2010 Mathematics Subject Classification: Primary: 52C07, Secondary: 52A40


\section{Preliminaries}

We will work in $n$-dimensional Euclidean space $\bR^n$. The norm
of a vector $x$ is denoted by $\|x\|$. For concepts not
defined in this paper we refer, e.g., to \cite{Gruber-Lekkerkerker} and
\cite{Rogers}.

A set $K \subset \bR^n$ is a \emph{convex body} if it is convex,
compact, and its \emph{interior} int\,$K$ is non-empty. For $A \subset \bR^n$ we
denote by $V(A)$ its \emph{volume} (Lebesgue measure; supposing it
exists -- but we will use this for bounded convex sets only), and
by lin$\,A$, aff$\,A$, or conv$\,A$ its \emph{linear}, \emph{affine},
or \emph{convex hull},
respectively. The \emph{dimension} dim$\,A$ of a convex set $A \subset \bR^n$
is the dimension of its affine hull.
$B^n \subset \bR^n$ is the unit
ball of $\bR^n$ with centre $0$. Its volume $V (B^n)$ is denoted by
$\kappa_n$. We
have $\kappa_n = \pi^{n/2}/\Gamma(1+n/2) \sim n^{-n/2} (2 \pi
e)^{n/2} \sqrt{\pi n}$ (by  Stirling's formula). If $X \subset
\bR^n$ is a linear subspace, $X^\perp \subset \bR^n$ denotes its
\emph{orthocomplement}. For $x,y \subset \bR^n$ we denote by
$\langle x,y \rangle$ their \emph{scalar product}. If $K \subset
\bR^n$ is a convex body with $0 \in {\rm int} K$, then its
\emph{polar} (with respect to the unit sphere about $0$) $K^*$ is
defined as $\{y \in \bR^n \mid \forall x \in K \, \langle x,y \rangle \le
1\}$, see \cite{Gruber-Lekkerkerker}, p. 107. In our paper we will deal
with \emph{lattices} in $\bR^n$, which will usually be denoted by
$L$ (and sometimes by $\Lambda$).
A {\emph{cross-polytope}} in $\bR^n$ is an affine image of the set conv$\{ \pm
e_1,..., \pm e_n \} $, where $e_1,...,e_n$ is the standard base of $\bR^n$.

There is a natural topology on the set of all lattices in $\bR^n$.
\emph{A neighbourhood base of a lattice $L \subset \bR^n$} is
obtained in the following way. The sets $\{\Lambda \subset \bR^n \mid
\Lambda$ is a lattice in $\bR^n$ and has a base $\{y_{1}, \dots, y_n\}$
with $\|y_1 - x_1 \| < \varepsilon, \dots, \| y_n - x_n \| <
\varepsilon\}$, where $\{x_1, \dots, x_n\}$ is any (some) base of $L$,
and $\varepsilon \in (0,\infty)$, is a neighbourhood base of $L$.
The definitions with ``any'' and ``some'' are equivalent; in
particular, the definition with ``some'' is independent of the
choice of the base of $L$ in question. This topology is metric and
locally compact. The set of all those lattices whose minima are at
least $c\,\,(>0)$ and the absolute  values of whose determinants are
at most some $C\,\,(<\infty)$ is compact (for all these facts, cf. 
\cite{Gruber-Lekkerkerker}, pp. 177-180).


\begin{definition}\label{def1.1}
Let $0 \le k \le n-1$ be integers, and $K \subset \bR^n$ be a convex
body. A lattice of translates of $K$ is $k$-\emph{impassable} if each
affine $k$-subspace of $\bR^n$ meets some body of the body
lattice. (For $k = 0$ this is the well known concept of lattice
covering. For $k=1$ or $k=n-1$, $1$-impassable or
$(n-1)$-impassable are also called \emph{non-transilluminable} or
\emph{non-separable}, respectively.)
\end{definition}


R. Kannan and L. Lov\'{a}sz \cite{Lovasz-Kannan} also
investigate this property and introduce, connected to this, covering
minima of a convex body with respect to a lattice. However, G.
Fejes Tóth's question, to which we refer here, and the investigations of
R. Kannan and L. Lov\'{a}sz 
\cite{Lovasz-Kannan} seem to go into practically disjoint
directions, except
for Lemmas (1.2) and (2.3) in \cite{Lovasz-Kannan}.

The main subject of this paper will be the investigation of a question
of G. Fejes Tóth; see \cite{Toth-1997}. This is the following. Let
$0 \le k \le n-1$ be integers. Find the infimum of the densities of
lattices of closed unit balls that are $k$-impassable.
We will also investigate the variant when $B^n$ is replaced by a
fixed convex body $K \subset \bR^n$. Rather up-to-date results and
problems about this concept are discussed in \cite{Brass-Moser-Pach}, pp.
149-159.


\begin{definition}\label{def1.2}
Let $0 \le k \le n-1$ be integers. For $K \subset \bR^n$ a
convex body, $d_{n,k}(K)$ will denote the \emph{infimum of
densities of $k$-impassable lattices of translates of $K$}. If $K =
B^n$, we write $d_{n,k} (B^n) = d_{n,k}$.
\end{definition}


Evidently, $d_{n,k} (K)$, as a function of $K$, is affine
invariant. Also, we have evidently $d_{n,0} (K) \ge \dots \ge
d_{n,n-1} (K)$. Since $d_{n,0} (K)$ is the well known quantity called
the density of the thinnest lattice covering by translates of $K$, we will
investigate in general the case $1 \le k \le n-1$ only.

Let $L \subset \bR^n$ be a lattice. By $D(L)$ we denote \emph{the
absolute value of its determinant} or, what is the same, the
\emph{volume of a basic parallelotope of $L$}.

J. M. Wills \cite{Wills} introduced a generalization of this
concept. For $k=1$ this concept is the well-known concept of the
minimal length of a (non-$0$) lattice vector.


\begin{definition}\label{def1.3}
(see \cite{Wills})
Let $1 \le k \le n-1$ be integers. For $L \subset \bR^n$ a lattice,
$D_k (L)$ will denote the minimum of (the absolute values of) the
determinants of its $k$-dimensional sublattices. (Of course, we may
suppose that we consider only sublattices of the form $L \cap
X_k$, where $X_k$ is a linear subspace spanned by some $k$ vectors
of $L$.)
\end{definition}


\begin{definition}\label{def1.4}
(H. Minkowski, e.g. \cite{Gruber-Lekkerkerker}, p. 58)
Let $1 \le k \le n$ be integers, $K \subset \bR^n$ a $0$-symmetric convex body,
and $L \subset \bR^n$ a lattice. The \emph{successive minima of $K$ with
respect to} $L$,
denoted by $\lambda_k (K,L)$, are defined by
\begin{equation}\label{equa1.1}
\lambda_k (K,L) = \min \{\lambda > 0 \mid \dim  {\rm lin} (\lambda K
\cap L) \ge k\}\,.
\end{equation}
\end{definition}

We will use this concept only for $K = B^n$.


\begin{definition}\label{def1.5}
(Dirichlet-Voronoi)
Let $A \subset \bR^n$ be such that inf $\{\|a_1 - a_2\| \mid a_1
\not= a_2 \in A\} > 0$ and there exists a positive number so that
any closed ball in $\bR^n$ with that radius intersects $A$. Then
the \emph{Dirichlet-Voronoi cell (D-V cell) of $a \in A$ with respect to} $A$
is defined as $\{x \in \bR^n \mid \forall b \in A \setminus \{a\}
\,\,\,
\|x-a\| \le \| x-b\|\}$. This is a convex polytope in $\bR^n$.
\end{definition}
We will use this concept only for $A$ a lattice.


\begin{definition}\label{def1.6}
Let $n \ge 1$ be an integer, and $K \subset \bR^n$ be a convex body.
We denote by $\delta_L (K)$, or $\vartheta_L (K)$, \emph{the density of the
densest lattice packing in $\bR^n$, or of the thinnest lattice covering of 
$\bR^n$,
by
translates of} $K$, respectively.
\end{definition}


For $K = B^n$, the density
$\delta_L (B^n)$ is known for $k \le 8$ (see \cite{Toth-1997}, p. 23), 
and recently
it has been announced for $k = 24$ by H. Cohn and A. Kumar \cite{Cohn-Kumar1}, 
with an outline of proof,
and proved by the same authors in \cite{Cohn-Kumar2},
while $\vartheta_L (B^n)$ is known for $n \le 5$ (see
\cite{Toth-1997}, p. 23). These are
well-investigated quantities, and we will treat them as ``known
quantities''.

The following estimates are known for $\delta_L (K), \vartheta_L
(K), \delta_L (B^n), \vartheta_L (B^n)$, see \cite{Toth-1997}, pp.
149-150:
\begin{equation}\label{equa1.2}
\frac{(\log 2-\varepsilon) \sqrt{\pi} n^{3/2}}{4^n} \le
\delta_L (K) \le 1 \le \vartheta_L (K) \le n^{(\log \log n)/\log 2 +
{\rm const}}\,,
\end{equation}
the left hand side inequality holding for $n \ge
n_\varepsilon$ (with $n_\varepsilon$ an integer depending on
$\varepsilon$). For centrally symmetric $K$ the first inequality
can be sharpened to
\begin{equation}\label{equa1.3}
\frac{(\log 2-\varepsilon)n}{2^n} \le \delta_L (K)\,,
\end{equation}
for $n \ge n_\varepsilon$.

For $B^n$ one has better estimates (\cite{Toth-1997} p. 23,
\cite{Toth-1999}, pp. 149-150, \cite{Gruber-Lekkerkerker},
p. 505, and \cite{Rogers}, p. 19):
\begin{equation}\label{equa1.4}
\begin{array}{lll}
2 \zeta (n) (n-1)/2^n \le \delta_L (B^n) \le
2^{-(0,5990... + o(1))n}\,,
\end{array}
\end{equation}
\begin{equation}\label{equa1.5}
\begin{array}{l}
n/(e^{3/2}+o(1)) \le \vartheta_L (B^n) \le {\rm const} \cdot n \cdot
(\log n)^{\left( \log (2 \pi e) \right) /\log 4}\,.
\end{array}
\end{equation}

Our paper has a non-trivial overlap with lower density estimates of
$k$-impassable lattices of convex bodies in \cite{MH} (in particular, their
Theorem 3.1), but the results of
their and our paper have been obtained independently.


\section{Introduction}

In this section we will recall earlier results on the subject of
our paper. L. Fejes T\'oth and E. Makai, Jr.\cite{Toth-Makai} stated
\begin{equation}\label{equa2.1}
d_{2,1} = \sqrt{3} \pi/8\,,
\end{equation}
\begin{equation}\label{equa2.2}
\min \{d_{2,1} (K) \mid K \subset \bR^2 \, \mbox{ is a convex body}\} =
3/8\,.
\end{equation}
In fact, they  proved only that the left hand sides
are at least the right hand sides, but did not prove that in the claimed cases
of
equalities in fact equalities hold, thus that the estimates are sharp.
However, these follow from E. Makai, Jr.
\cite{Makai-1978}. Namely, \cite{Makai-1978}
proved
\begin{equation}\label{equa2.3}
\begin{array}{l}
\min \{ d_{2,1} (K) \mid K \subset \bR^2 \mbox{ is a centrally
symmetric convex body} \}  \\
= 1/2\,.
\end{array}
\end{equation}
In fact, moreover \cite{Makai-1978} 
claimed that the minimum is attained only for
a parallelogram, for which he used that
\begin{equation}\label{equa2.4}
\min \{V(K) V (K^*) \mid K \subset \bR^2 \, \mbox{ is
0-symmetric convex body}\} \, = 8
\end{equation}
(a theorem of K. Mahler \cite{Mahler39}), the minimum being attained only for a
parallelogram. This last statement has been proved by
\cite{Reisner}, thus completing the proof of the theorem in
\cite{Makai-1978}.

We have to mention that a sketch of the proof of (\ref{equa2.2}) was
also given in \cite{Makai-1978}, see the proof of Theorem 3. For a proof of
the inequality claimed and used in the proof of (\ref{equa2.2}) in
\cite{Makai-1978}, proof of Theorem 3, we refer to \cite{Eggleston}.

R. Kannan and L. Lov\'asz
\cite{Lovasz-Kannan}, Lemma (2.3), and \cite{Makai-1978}, Theorem 1 
proved for $K \subset
\bR^n$ a convex body
\begin{equation}\label{equa2.5}
d_{n,n-1} (K) = \frac{V(K) V \left( \left( (K-K)/2 \right) ^* \right) }
{4^n \delta_L \left( (K-K)/2 \right) }\,.
\end{equation}
In particular,
\begin{equation}\label{equa2.6}
d_{n,n-1} = \frac{\kappa^2_n}{4^n \delta_L (B^n)}\,.
\end{equation}
Since $\delta_L (B^n)$ is known for $n \le 8$ (cf. \cite{Toth-1997}, p. 23),
and for $n = 24$ has been
proved in \cite{Cohn-Kumar2}, $d_{n,n-1}
(B^n)$ is known for $n \le 8$ and 
for $n=24$. (Here we have to make a correction to
\cite{Makai-1978}.
There it was stated that for $n=3$ the thinnest non-separable lattice of
translates of $B^3$ is given by a lattice $L$ similar to the lattice or
thinnest lattice
ball covering, i.e., a space-centered cubic lattice. This holds true, but the
ratio of
similarity was given incorrectly. $L$ has correctly a base $\{\sqrt{2}
(-1,1,1), \sqrt{2} (1,-1,1),
\sqrt{2} (1,1,-1)\}$.)
As a consequence, cf. \cite{Makai-1978}
for $K \subset \bR^n$ a convex body, or an $0$-symmetric convex body,
\begin{equation}\label{equa2.7}
d_{n,n-1} (K) \ge V (K) V (((K-K)/2)^*) 4^{-n}\, \mbox{ or }
\end{equation}
\begin{equation}\label{equa2.8}
d_{n,n-1} (K)\ge V(K) V (K^*)4^{-n}\,,
\end{equation}
respectively.

\cite{Lovasz-Kannan} proved
\begin{equation}\label{equa2.9}
\begin{array}{lll}
\min \{d_{n,n-1} (K) \mid K \subset \bR^n \, \mbox {is a centrally
symmetric convex body}\} \ge \\
\min \{ V(K)V(K^*)4^{-n}
\mid K \subset \bR^n \, \mbox {is a }\,\,0 \mbox {-symmetric convex body}\}
\\
=n^{-n} e^{O (n)}\,.
\end{array}
\end{equation}
Here for the last inequality they used the theorem of
J. Bourgain-V. D. Milman, cf.
\cite{Bourgain-Milman}. However, this theorem was recently improved by
G. Kuperberg \cite{Kuperberg}, namely to
\begin{equation}\label{equa2.9A}
V(K)V(K^*) > \kappa _n^2/2^n\,.
\end{equation}
So, combining (\ref{equa2.9}) and (\ref{equa2.9A}), we have
\begin{equation}\label{equa2.9B}
\begin{array}{l}
\min \{ d_{n,n-1} (K) \mid K \subset \bR^n \mbox{ is a centrally
symmetric convex body} \} \ge\\
\min \{  V(K)V(K^*)4^{-n}
\mid K \subset \bR^n \mbox{ is a }\, 0 \mbox{-symmetric convex body} \} 
\\
> \kappa _n^2/8^n\,.
\end{array}
\end{equation}

However, from (\ref{equa2.7}) and (\ref{equa2.9A}),
together with the observation that, for $K
\subset \bR^n$ a convex body, $V(K) V(((K-K)/2)^*) = V(K) V((K-K)/2)^{-1}
\cdot V((K-K)/2) V(((K-K)/2)^*)$, where
the first factor is at least $2^n {2n \choose n}^{-1}$
(see \cite{Rogers}, p. 37), and at most $1$ (by the
Brunn-Minkowski inequality, cf. \cite{Schneider}, {\S} 6.1)
we have
\begin{equation}\label{equa2.10}
\begin{array}{l}
\min\{d_{n,n-1} (K) \mid K \subset \bR^n \, \mbox{ is a convex body}
\} \ge\\
\min\{V(K) V(((K-K)/2)^*)4^{-n}\mid K \subset \bR^n \, \mbox{ is a convex
body} >\\
\kappa _n^2  {2n \choose n}^{-1}4^{-n}\,.
\end{array}
\end{equation}

There is a
counterpart to the question for the densest lattice packing of a
convex body $K \subset \bR^n$ (that is, of determining $\delta_L
(K)$). Namely, one can ask for an inequality from the other side:
find $\min \,\{\delta_L (K) \mid K \subset \bR^n$ is a (centrally
symmetric) convex body$\}$. These theorems are the so-called Minkowski-Hlawka
type theorems; cf. \cite{Gruber-Lekkerkerker}, \S 19.5 Theorem 8,
for general $n$ and general convex bodies, and
\S 22, Theorem 7, for $n=2$ and general convex bodies. 
We will need only the case $n=2$,
centrally symmetric bodies. For them we have
the inequality of P. Tammela \cite{Tammela}:
\begin{equation}\label{Tammela}
\begin{array}{l} 
\min \,\{\delta_L (K) \mid K \subset \bR^2 {\mbox{ is a 
centrally symmetric convex body}} \} \\
\ge 0.8926... \,.
\end{array}
\end{equation}
\cite{Tammela} has at its right hand side a certain,
explicitly given algebraic number, which is however defined in a quite
complicated way.
Therefore we do not reproduce it here. However, since \cite{Tammela}
is in Russian, we give
hints to this definition. The last displayed formula in \cite{Tammela} is the
inequality cited by us.
The quantity $m$ in it is defined by \cite{Tammela}, (18), (22), as a certain
rational function of $u_0$, where
$u_0$ is the unique root of the polynomial equation (23) in \cite{Tammela}, in
the interval given in
(23) from \cite{Tammela} (the total length of these formulas is twelve lines).

Analogously, for $d_{n,k}(K)$ there is also the question 
of opposite, i.e., upper estimates, i.e., to a minimax problem.
\cite{Makai-1978}, Theorems 4 and 5 proved that
\begin{equation}\label{equa2.11}
\begin{array}{l}
\max\{d_{2,1}(K) \mid K \subset \bR^2 \, \mbox{ is a convex body}\} =
\\
\max\{d_{2,1}(K) \mid K \subset \bR^2 \, \mbox{ is a centrally
symmetrix convex body}\}  \\
\le \pi^2 / [4(3 \sqrt{2} + \sqrt{3} -
\sqrt{6})] = 0,6999\dots\,,
\end{array}
\end{equation}
and
\begin{equation}\label{equa2.12}
\begin{array}{l}
\max \{d_{n,n-1} (K) \mid K \subset \bR^n  \mbox{ is convex body}\}
=\\
\max\{d_{n,n-1} (K) \mid K \subset \bR^n  \mbox{ is a centrally
symmetric convex body}\}\\
 \le \kappa ^2_n /[(n \log 2 - {\rm const)}
2^{n-1}]\,,
\end{array}
\end{equation}
for $n$ sufficiently large.

Concerning $d_{n,k}$, besides the case $k = n-1$ there is just one case proved
by R. P. Bambah, A. C. Woods
\cite{Bambah-Woods}, namely:
\begin{equation}\label{equa2.13}
d_{3,1} = 9 \pi/32\,.
\end{equation}
We note that the minimum is attained, e.g., for a lattice similar
to the lattice of the densest lattice packing of balls in $\bR^3$.
(\cite{Bambah-Woods} contains a small gap. In p. 153 it considers a quadratic
form $G(x_1,x_2)=ax_1^2+2bx_1x_2+cx_2^2$, assuming integer values for all
integers
$x_1,x_2$. It states that $a,b,c$ are integers. However, in fact, only
$a,2b,c$ are integers. This affects case (i) in p. 153 and case (ii) in
p. 154. In case (ii) we have the additional cases $a=4$, $b=1/2$, and $a=4$,
$b=3/2$, both of which lead to a contradiction like in \cite{Bambah-Woods}.
In case (i) we have the additional cases $a=3$, $b=1/2$, and $a=3$,
$b=3/2$. Proceeding like in \cite{Bambah-Woods}, we do not get a
contradiction, but we obtain, as unique possibility, $c=3$ in both of these
cases. Then the determinant of (the symmetric matrix associated to)
our quadratic form $G$ is $35/4$, or $27/4$. The
norm square of the primitive vector $OA_1$, in p. 138, is the quotient of the
determinants of the
quadratic forms associated to the investigated $3$-dimensional lattice,
and of its $2$-dimensional projection. That is, with the notations of 
\cite{Bambah-Woods}, $\| OA_1 \| ^2= a_{1,1}d(f)/d(G)$, where $a_{1,1}=2$ and
$d(f)=4$, i.e., $\| OA_1 \| ^2$ is $64/35$,
or $27/4$. However, the norm squares of the vectors of the
investigated $3$-dimensional lattice are even integers, a contradiction.)

Next we turn to a subject which is a bit different, when we will consider,
rather than any lattices of translates of a convex body, only
lattice packings of translates of a convex body. Here we will list
some results for the unit ball. The relation of these two types of
questions will become clear later in this paper.

A. Heppes \cite{Heppes} proved the following:
\begin{equation}\label{equa2.14}
\begin{array}{l}
\mbox{For each lattice packing of closed unit balls in } \, \bR^3 \mbox{ there
exist}\\
\mbox{three both-way infinite open circular cylinders, with
linearly}\\
\mbox{independent
directions of axes, disjoint to the union of the}\\
\mbox{ball lattice}.
\end{array}
\end{equation}

In fact, the proof of \cite{Heppes} shows that the radii of the
bases can be $2 /\sqrt{3}-1$. I. Hortob\'agyi \cite{Hortobagyi} stated the
following sharpening of the above result:
\begin{equation}\label{equa2.15}
\begin{array}{l}
\mbox{For each lattice packing of closed unit balls in }\, \bR^3 \mbox{ there
exist}\\
\mbox{three both-way
infinite open circular cylinders with linearly}\\
\mbox{independent directions of
axes, disjoint to the union of the ball} \\
\mbox{lattice, with radii of their bases } \, 3 \sqrt{2}/4-1 = 0,0606\dots\,.
\mbox{ This}\\
\mbox{inequality is sharp, even when stated for one cylinder only,}\\
\mbox{for the densest lattice packing of closed unit balls}\,.
\end{array}
\end{equation}
\cite{Hortobagyi} did not prove that in the claimed case of
equality in fact equality holds (even when stated for one cylinder only),
thus that his estimate is
sharp. We will see that the result of \cite{Bambah-Woods} and our
Theorem 3.25 together will prove that in the above 
claimed case of equality in fact equality holds; not even one such
cylinder exists with a greater radius of its base.

J. Horv\'ath \cite{Horvath-1970}, Satz 1, and J. Horv\'ath and S. S. Ry\v skov
\cite{Horvath-Ryskov-1975} proved an analogous
result for $\bR^4$:
\begin{equation}\label{equa2.16}
\begin{array}{l}
\mbox{For each lattice packing of closed unit balls in } \bR^4 \mbox{ there
exist}\\
\mbox{four both-way infinite open cylinders, disjoint to the union
of}\\
\mbox{the ball lattice,
with linearly independent directions of axes}\\
\mbox{and with bases that are 3-balls of radii }
\sqrt{5}/2-1 \ge 0,1180\dots\,.\\
\end{array}
\end{equation}
(We note that still in \cite{Horvath-1970} sharpness of  the result
was claimed, and this was repeated in \cite{Ryskov-Horvath}, p. 128;
however, this was withdrawn in
\cite{Horvath-Ryskov-1975},
available only in Hungarian, in p. 92, (2), where already only
``$R_4 \ge \sqrt{5}/2-1$'' stands, together with an explanation
in p. 94, Paragraph 2, end of last line, that ``$R_4$ is probably
sharp''.)

A critic of \cite{Horvath-1970}, Satz 2, is contained in T. Hausel
\cite{Hausel}. He shows:
\begin{equation}\label{equa2.17}
\begin{array}{l}
\mbox{There exist lattice packings of closed balls in} \, \bR^n
\,\mbox{such that for}\\
k \ge n- c \sqrt{n} \,\, (c > 0 \mbox{ some constant) they are
$k$-impassable}\,.
\end{array}
\end{equation}
This was sharpened by \cite{Henk-Ziegler-Zowa}, or \cite{Henk},
who showed that
\begin{equation}\label{equa2.18}
\begin{array}{lll}
\mbox{in the above statement one may replace the hypothesis} \\
\,k \ge n-c\sqrt{n} \mbox{ by } k \ge n - cn \,\,(c > 0
\mbox{ some constant)}, \\
\mbox{or by } n \mbox{ is sufficiently large, and } k \ge n/\log_2 n,\\
\mbox{and still the same conclusion holds}\,.
\end{array}
\end{equation}


\section{Results}

 First we determine $d_{n,n-1}(K)$ for a
simplex. (This was claimed, but not proved, in \cite{Makai-1978},
Proposition 1. The statements there, not included in our Proposition 1
here, easily follow from our proof here.)


{\textbf{Proposition 3.1.}}
\emph{Let $n \ge 2$, and let $K \subset {\Bbb R}^n$ be an $n$-simplex. 
Then we have}
\[
d_{n,n-1} (K) = \frac{n+1}{2^n n!} = n^n
\left(\frac{e}{2}\right)^n e^{o(n)}\,.
\]
\vskip.2cm


{\textbf{Definition 3.2.}}
(\cite{Wills})
Let $1 \le k \le n-1$ be integers, and $L \subset \bR^n$ be a
lattice. We denote by $D_k(L)$ the quantity $\min\{|\det L_k| \mid
L_k$ is an $k$-dimensional sublattice of $L\}$. (Evidently we may
assume additionally that $L_i = L \cap {\rm lin}\,L_i$; this does not affect
the value of the minimum.) For $k = d$ we write $D(L) = D_d(L)$. For $d=1$, we
have that
$D_1(L)$ is the minimal length of a non-$0$ vector in $L$.
\vskip.2cm


The existence of the minimum
was stated both in \cite{Wills}, p. 268, and \cite{Schnel}, p. 607.
It follows from the version of the selection theorem of Mahler
(\cite{Gruber-Lekkerkerker}, p. 179) for $k$-dimensional lattices
in $\bR^n$, for which the proof in \cite{Gruber-Lekkerkerker}, pp.
179-180, goes through without any change. Then the hypotheses of
this version of Mahler's Theorem are satisfied. 1) Each
$k$-dimensional sublattice of $L$ has as minimum at least that of $L$.
2) Looking for the infimum of
the absolute values of the $k$-dimensional determinants of the
$k$-dimensional sublattices, we may assume that we consider only such
$k$-dimensional sublattices for which these absolute values of
$k$-dimensional determinants are bounded from above.


{\textbf{Definition 3.3.}}
Let $1 \le k \le n-1$ be integers. We denote by $c_{n,k}$ the
number $\min \{c \in (0,\infty) \mid $ for each lattice $L \subset
\bR^n$ we have $D_k (L) \le c \cdot D(L)^{k/n}\}$.
\vskip.2cm


In Theorem 3.5 we will see that the set of numbers $c$ in
question is nonempty. Then the existence of the minimum of the
numbers $c$ is obvious.

We remark that $c_{n,1}$ is a well-known quantity: we have
\[
\left(\frac{c_{n,1}D(L)^{1/n}}{2}\right)^n \kappa_n / D(L) = \delta_L
(B^n)\,.
\]
From here and from the estimates in \S~1 we have
\[
\sqrt{n} \left(\sqrt{\frac{2}{\pi e}} \cdot \frac{1}{2} - o(1)\right)
\le c_{n,1} \le \sqrt{n} \left(\sqrt{\frac{2}{\pi e}} \cdot
\frac{1}{2^{0,5990\dots}}+ o(1)\right)\,.
\]

On the other hand, \cite{Schnel}, Theorem 3, says
\begin{equation}\label{equa3.1}
c_{n,n-1} \le (n^{3/2}/2)\kappa_n^{1/n} \sim \sqrt{\pi e/2} \cdot n\,,
\end{equation}
which is much better than our following theorem, applied to the case $k
= n-1$. (Also, (\ref{equa3.1}) substantially improves \cite{Wills},
Theorem 3.) Our Theorem 3.5 will be good for not too large values
of $k$, cf. Remark 3.6. If we would have defined also $c_{n,n}$,
analogously we would have $c_{n,n}=1$. Anyway, for each $k \in
\{1,\dots,n-1\}$ we have $c_{n,k} \ge 1$. In fact, for $L={\bZ}^n$ we
have $D_k ({\bZ}^n)$ a positive integer. And, since ${\bZ}^k \subset {\bZ}^n$,
actually $D_k ({\bZ}^n)=1$. Observe still that in the next theorem the
$(1/k)$-th power of the upper estimate for $c_{n,k}$ is independent
of $k$.


{\textbf{Problem 3.4.}}
Determine as many $c_{n,k}$'s as possible.
\vskip.2cm


{\textbf{Theorem 3.5.}}
\emph{Let $1 \le k \le n-1$ be integers. Then the numbers $c_{n,k}$ from
Definition 3.2 satisfy
\[
c_{n,k} \le 2^k \left(\frac{\delta_L(B^n)}{\kappa_n}\right)^{k/n} \le
2^k \left(2^{-0,5990\dots} \cdot \sqrt{\frac{n}{2\pi e}}
(1+o(1))\right)^k\,.
\]
For $k = 1$ the first inequality is an equality, while for $2
\le k \le n-1$ it is a strict inequality.}
\vskip.2cm


The question whether $d_{n,k}$, or $d_{n,k} (K)$, respectively, is
a minimum is not clear. For $k=0$ we have the question of
existence of the thinnest lattice covering with translates of $B^n$,
or $K$, respectively, which is known to exist. For $k=n-1$, by
\cite{Lovasz-Kannan} and \cite{Makai-1978} a lattice $\{x+K\mid x\in
L\}$ of translates of a convex body $K$ is non-separable if and only if the
body lattice $\{x^* +
(((K-K)/2)^*)/4 \mid x^* \in L^*\}$,
where $L^*$ is the lattice polar to $L$,
is a lattice packing, and the
product of the densities of these two body lattices is independent
of $L$ (it depends on $K$ only). Since densest lattice packings of a
convex body exist, there exist also thinnest non-separable
lattices of a convex body.

However, for $1 \le k \le n-2$ the situation seems to be
different. We have


{\textbf{Proposition 3.6.}}
\emph{Let $1 \le k \le n-2$ be an integer, $K \subset \bR^n$ a strictly
convex body, and $X \subset \bR^n$ a $k$-dimensional linear
subspace. Then the set of lattices $L \subset \bR^n$ such that for the body
lattice
$\{x + K\mid x \in L\}$ there
does not exist a translate of $X$ disjoint to $\bigcup\{x+K\mid x\in L\}$
is not closed in the topology of lattices described in \S~1.
Moreover, there exists a lattice $L$ such that there is a
translate of $X$ disjoint to $\bigcup\{x+K\mid x \in L\}$, but
for any $\varepsilon >0$ there does not exist any translate of $X$ disjoint
to $\bigcup\{x+K\mid x \in (1-\varepsilon)L\}$
(and here $(1-\varepsilon) L$ converges to $L$ for $\varepsilon \to
0$)}.
\vskip.2cm


For the quantity $d_{n,k}$, defined in Definition 1.2, we
have


{\textbf{Theorem 3.7.}}
\emph{Let $1 \le k \le n-1$ be integers. Then we have
\[
\begin{array}{c}
d_{n,k} \ge 
\frac{\kappa_n \vartheta_L (B^{n-k})^{n/(n-k)}}{
\kappa_{n-k}^{n/(n-k)}c_{n,k}^{n/(n-k)}}
\ge \\
\frac{\kappa_n \vartheta_L (B^{n-k})^{n/(n-k)}} 
{\kappa_{n-k}^{n/(n-k)}
2^{kn/(n-k)}
\left( \delta_L(B^n)/\kappa_n \right)^{k/(n-k)}}\,.
\end{array}
\]
For $k \ge 2$ the second inequality is strict. 
For $k=1$, in the inequality, obtained from the above chain of inequalities,
by omitting the middle term, we have equality for $n=2$, and we have strict
inequality for $3 \le n \le 6$. 
We have also
\[
d_{n,k} \ge \frac{\kappa ^2_n}{4^n \delta_L (B^n)} = \frac{e^{O(n)}}{n^n}\,.
\]}
\vskip.2cm


{\bf{Remark 3.8.}}
The smaller lower estimate in the first chain of inequalities in
Theorem 3.7 can be still (not substantially) diminished by using
$\vartheta_L (B^{n-k}) \ge (n-k)e^{-3/2}(1+o(1))$, and $\delta_L
(B^n) \le 2^{-(0,5990 \dots + o(1))n}$, and then by Stirling's
formula we obtain, for $n \to \infty $ and $n-k \to \infty$,
\[
\begin{array}{c}
d_{n,k} \ge \left( \left( 1/(\sqrt{n})^{k/(n-k)} \right)
\left(\sqrt{\pi e/2} \cdot 2^{0,5990\dots + o(1)} \right) ^{k/(n-k)}
\sqrt{(n-k)/n} \right) ^n  \\
\times \left( (n-k)e^{-3/2} \right) ^{n/(n-k)}
\left( (n-k)/n \right) ^{n/(n-k)}\,.
\end{array}
\]
Here the factors are written in the order according to their
``contribution to the product''. Thus one sees, e.g., that for $k
\ge 1$ fixed, and $n \to \infty$, this behaves approximately like
${\rm const}_k \cdot n / (\sqrt{n})^k$. If $k/n=c < 2/3$ is
fixed, it behaves like $e^{O(n)}/n^{nc/(2-2c)}$, where
$c/(2-2c)<1$, which is still a better estimate than
$\kappa^2_n/(4^n\delta_L(B^n))=n^{-n} e^{O(n)}$. We yet
remark that, for $k \ge n-1$, the first inequality of
Theorem 3.7 and (\ref{equa3.1}) give $d_{n,n-1} \ge
n^{-3n/2}e^{O(n)}$, rather than the correct $d_{n,n-1} =
n^{-n} e^{O(n)}$, cf. (\ref{equa2.6}).
\vskip.2cm

Since $\delta_L(B^n)$ is known
for $n \le 8$ and $n = 24$, and $\vartheta_L (B^n)$ is
known for $n \le 5$, we can evaluate the second lower bound in
Theorem 3.7 for $n \le 8$ and $n=24$, and
for $n-k \le 5$. Since the second lower bound in
Theorem 3.7 is sharp only for $k=1$, we can have sharp results from
Theorem 3.7 only for $k=1$, and thus $n \le 6$.

Observe that for the case $n = 2$ the sharp estimate
$d_{2,1} = \sqrt{3} \pi/8$ was given already in \cite{Toth-Makai}, also cf.
\cite{Makai-1978}. 
As concerns the case $k=1$ and $n \ge 3$, in
\cite{Horvath-Ryskov-1975}, p. 92, (2), it was stated, without proof, that for
$n = 5,6$ the projection of the (unique; see
\cite{Gruber-Lekkerkerker}, p. 517) densest lattice packings of
balls in $\bR^n$, when projected
along any minimum vector, will not yield a lattice similar to a
lattice of thinnest lattice coverings of balls in $\bR^{n-1}$. We will give a
proof of this for ``some minimum vector'' 
in the proof of Theorem 3.7, which will imply the
non-sharpness of the estimates mentioned in Theorem 3.7, for $n=5,6$. (For
$n=5$ the minimum vectors are equivalent under the group of congruences of 
the respective lattice.)
For $n = 3$ the analogous statement is simple; cf. in the proof of Theorem
3.7.
For $n=4$ we do have similarity of the two lattices (the
projection of a space centred cubic lattice along the direction of a
coordinate axis is a space centred cubic lattice, and all minimum vectors of
the space centred cubic lattice in ${\Bbb R}^4$ are equivalent under the group
of congruences of this lattice).
However, a suitable projection 
in a direction different from the
directions of all minimum vectors will prove the needed non-sharpness 
of the estimate mentioned in Theorem 3.7, for $n=4$.


We give an upper estimate for $d_{n,k}$, for certain $k$'s, which is probably
rather weak.


{\textbf{Proposition 3.11.}}
\emph{There exists a constant $n_0$ such that the following holds.
For any integers $n \ge n_0$ and $n/ \log _2 n \le k \le n-1$ we have
$d_{n,k} \le \delta_L (B^n) \le 2^{-(0,5990+o(1))n}$.}
\vskip.2cm


Now we turn from $B^n$ to arbitrary (centrally symmetric) convex
bodies. We use $c_{n,k}$ from Definition 3.3.


{\textbf{Theorem 3.12.}}
\emph{Let $1 \le k \le n-1$ be integers. Then we have
\[
\begin{array}{c}
\min\{d_{n,k} (K) \mid K \subset \,{\mathbb R}^n \mbox{ is a convex body} \}
\ge\\
\frac{d_{n,k}}{\kappa_n n!n^{n/2}/(n+1)^{(n+1)/2}} \ge \\
\frac{\kappa_n \vartheta_L (B^{n-k})^{n/(n-k)}/ \left( \kappa_{n-k}^{n/(n-k)}
c_{n,k}^{n/(n-k)}\right) } 
{\left[ \kappa_n n!n^{n/2} / (n+1)^{(n+1)/2}\right] } \ge \\
\frac{\kappa_n \vartheta_L (B^{n-k})^{n/(n-k)} / 
\left( \kappa^{n/(n-k)}_{n-k} 2^{kn/(n-k)}
\left( \delta_L (B^n)/\kappa_n \right) ^{k/(n-k)} \right) } 
{\kappa_n n!n^{n/2} / (n+1)^{(n+1)/2}}\,.
\end{array}
\]
Further, we have
\[
\begin{array}{c}
\min \{d_{n,k} (K) \mid K \subset {\mathbb R}^n \mbox{ is a centrally
symmetric convex body} \} \ge\\
\frac{d_{n,k}}{\kappa_n n!/2^n} \ge \\
\frac{\kappa_n \vartheta_L (B^{n-k})^{n/(n-k)} / (\kappa^{n/(n-k)}_{n-k}
c^{n/(n-k)}_{n,k})}{\kappa_n n! / 2^n} \ge \\
\frac{\kappa_n \vartheta_L (B^{n-k})^{n/(n-k)} / \left[ \kappa^{n/(n-k)}_{n-k}
2^{kn/(n-k)} \left(\delta_L (B^n)/\kappa_n\right)^{k/(n-k)} \right] } 
{\kappa_n n!/2^n}\,.
\end{array}
\]
Also we have
\[
\min \{ d_{n,k} (K) \mid K \subset \bR^n {\mbox{ is a convex
body}} \} \ge 
\frac{\kappa _n^2}{{2n \choose n}}=\frac{e^{O(n)}}{n^n}\,.
\]}
\vskip.2cm


{\textbf{Remark 3.13}}
The denominators of the fractions in Theorem 3.12 in
square brackets are $n^{n/2} e^{O(n)}$. Hence, for $k/n = c <
1/2$ fixed, the smaller lower estimate in the chain of
inequalities in Theorem 3.7 behaves ``approximately''
like $1/n^{c/(2-2c)}$, where $c/(2-2c) < 1/2$. Hence we get still
a better estimate than $n^{-n} e^{O(n)}$.
\vskip.2cm


Still we remark that for general $n,k$ we have no conjecture about
the infima of the quantities investigated in Theorem
3.12. However, for $k = n-1$ we do have (for the
existence of the minima, cf. the proof of Theorem 3.16):


{\textbf{Conjecture 3.14.}}
We have
\[
\min \{d_{n,n-1} (K)\mid K \subset \bR^n \, \mbox{ is a convex body}\}
\, = \frac{n+1}{2^n n!}\,,
\]
and
\[
\min \{d_{n,n-1} (K) \mid K \subset \bR^n \, \mbox{ is a centrally
symmetric convex body}\} \, = \frac{1}{n!}\,.
\]
If true, these would be sharp, cf. our Proposition 3.1
and \cite{Makai-1978}, Proposition 2, which show that in the first
case equality holds for a simplex, and in the second case for a
cross-polytope.
\vskip.2cm


Since by \cite{Makai-1978}, Theorem 1, we have $V(K) V ((K-K)/2)^*)
/ [4^n \delta_L (((K-K)/2)^*)] $
\newline
$\ge V(K) V((K-K)/2)^*]/4^n$, the
centrally symmetric case (when $((K-K)/2)^* = K^*)$ would follow
from the volume product conjecture of K. Mahler (i.e., that for
$0$-symmetric convex bodies $K \subset \bR^n$ we have $V(K) V(K^*)
\ge 4^n/n!)$; see \cite{Mahler39}.

Some particular cases of the centrally symmetric case can be proved.


{\textbf{Definition 3.15.}} (\cite{Goodey-Weil}, Goodey-Weil)
We call a convex body $K
\subset \bR^n$ a \emph{zonoid} if it is the Hausdorff limit of
some $K_i \subset \bR^n$, where each $K_i$ is a \emph{zonotope},
i.e., a finite vector sum of line segments. (Of course, each zonoid is
centrally symmetric.)
\vskip.2cm


For a $0$-symmetric convex body $K \subset
\bR^n$, the {\emph{associated norm}} is the norm on $\bR^n$
whose unit ball is $K$.


{\textbf{Theorem 3.16.}}
\emph{Let $n \ge 3$ be an integer. We have $\min\{d_{n,n-1}(K) \mid K
\subset \bR^n$ is a zonoid$\} > \min \{d_{n,n-1} (K) \mid K \subset
\bR^n$ is a polar of a zonoid centred at $0\} =
\min \{d_{n,n-1} (K) \mid K \subset
\bR^n$ is symmetric w.r.t. all coordinate hyperplanes$\} =
\min \{d_{n,n-1} (K) \mid K \subset
\bR^n$ is $0$-symmetric, and
the norm associated to $K$
satisfies that all the natural projections to the coordinate
hyperplanes are contractions$\} =
1/n!$.
For $n \le 8$ we have also $\min\{d_{n,n-1}(K) \mid K
\subset \bR^n$ is a $0$-symmetric convex polytope with at most
$2n+2$ facets$\} >
\min\{d_{n,n-1}(K) \mid K
\subset \bR^n$ is a $0$-symmetric convex polytope with at most
$2n+2$ vertices$ \}
= 1/n!$. The second, third and sixth minima
are attained only for cross-polytopes.}
\vskip.2cm


For the general case we have an
analogous


{\textbf{Conjecture 3.17.}}
For $K \subset \bR^n$ a convex body we have
\[
V(K) V(((K-K)/2)^*) \ge
\frac{2^n(n+1)}{n!}\,,
\]
possibly with equality only for a simplex.
(For a simplex we have equality, cf. the proof of Proposition
3.1.)
\vskip.2cm


The case $n=2$, including the case of equality 
(i.e., that it occurs only for
the triangle)
is proved by H. G. Eggleston \cite{Eggleston}. Observe that
\[
\begin{array}{c}
V(K) V \left( \left( (K-K)/2 \right) ^* \right) = \\
\left[ V(K)/V \left( (K-K)/2 \right) \right] 
\cdot \left[ V \left( (K-K)/2 \right) 
\cdot V \left( \left( (K-K)/2 \right) ^* \right) \right] \ge \\
\left[ 2^n \big/ {2n \choose n} \right] \cdot
\min \{ V(K)V(K^*) \mid K \subset \bR^n \mbox{ is a }
0 \mbox{-symmetric convex}\\
\mbox{body} \} \,.
\end{array}
\]
Here we used the difference body inequality; see \cite{Rogers}, p. 37.
If Mahler's volume product conjecture would be true, this quantity
would be $\sim[2^n(n+1)/n!] \sqrt{\pi/n}$, thus
very close to the conjectured value. Anyway, by (\ref{equa2.10}),
$V(K)V(((K-K)/2)^*) \ge \kappa _n^2
{2n \choose n}^{-1}$.
This remark hints to that the proof of
Conjecture 3.17 would be quite difficult (as is the case
with the sharp lower estimate in the volume product problem).

Analogously like above, in a particular case we have an almost sharp
estimate.


{\textbf{Theorem 3.18.}}
\emph{We have
$2^n (n+1)/n! \ge \min \{d_{n,n-1}(K) \mid K \subset \bR^n$ is a
convex body
such that $\left( (K-K)/2 \right) ^*$ is a zonoid$\} 
\ge [2^n / {2n \choose n}] / n! \sim [2^n(n+1)/n!] \sqrt{\pi/n}$.}
\vskip.2cm


Another way of approach for general convex bodies would be the following
by \cite{Makai-1978}, Corollary 1. For a lattice $L$ and a convex body $K$ the
body lattice $\{K+x\mid x \in L\}$ is $(n-1)$-impassable if and only if the
body lattice
$\{(K-K)/2+x\mid x \in L\}$ is $(n-1)$-impassable. Then include $(K-K)/2$ into
an ($0$-symmetric) ellipsoid
$E$ such that $V(E)/V(K)$ would be possibly small. Then also
$\{E+x\mid x \in L\}$ is $(n-1)$-impassable.
We have $d_{n,n-1} (K)/V(K) = d_{n,n-1} ((K-K)/2)/V((K-K)/2)$, and the ratio
of the densities of
$\{K+x\mid x \in L\}$ and $\{E+x\mid x \in L\}$ is $V(K) / V(E)$. Further, the
density of $\{E+x\mid x \in L\}$ is at
least $d_{n,n-1}$. So
\[
\begin{array}{c}
d_{n,n-1}(K) \ge\\
d_{n,n-1} /\min \{V(E)/V(K))\mid E
\subset \bR^n \mbox{ is an ($0$-symmetric) ellipsoid,}\\
(K-K)/2 \subset E\}\,.
\end{array}
\]


\textbf{Conjecture 3.19.}
For each convex body
$K \subset \bR^n$ there exists an ellipsoid $E \subset \bR^n$ such
that $(K-K)/2 \subset E$, and $V(E) / V(K)$ is at most the same quantity
when $K$ is a regular simplex and $E$ is the circumball of
$(K-K)/2$, i.e., a ball of the same diameter as $K$. In other
words, $V(E)/V(K) \le \kappa_n n!/(2^{n/2}\sqrt{n+1})$.


The case $n =
2$ is proved by \cite{Behrend}, p. 716, (II)$_3$, p. 715, (I), case
$\nu=5$.

This conjecture is very similar to the theorems of \cite{Ball},
\cite{Barthe-1997},
\cite{Barthe-1998}, and \cite{Barthe-2003}
used in the proof of
Theorem 3.12, but unfortunately does not follow from them. However, these imply
the following. By
\cite{Rogers}, p. 37, we have $V(K-K)/2)/ V(K) $
\newline
$\le \left(\begin{array}{c}2n\\n\end{array}\right) 2^{-n}$. By
\cite{Ball}, \cite{Barthe-1997}, \cite{Barthe-1998}, and
\cite{Barthe-2003}, cited at the beginning of the proof of Theorem
3.12, there exists an ellipsoid $E \supset (K-K)/2$ with
$V(E) / V((K-K)/2) \le \kappa_n n!/2^n$. Multiplying the two
inequalities, we have $V(E) / V(K) \le
\left(\begin{array}{c}2n\\n\end{array}\right) \kappa _n n!/4^n$. The
quotient of this upper bound and the conjectured one is $\sim
\sqrt{2^n/\pi}$.

Another way to formulate this conjecture is the following. To
look for an ellipsoid $E \supset (K-K) / 2$ means to look for an
affine image $K'$ of $K$ with $(K'-K')/2$ contained in $B^n$. Then
the quantity to be minimized is $\min \{V(B^n)/V(K')\mid K'$ is an
affine image of $K$ such that $(K'-K')/2 \subset B^n\}$. Observe
that this is a kind of ``reverse isodiametric inequality'' (cf. the
reverse isoperimetric inequality of \cite{Ball}). The isodiametric
inequality states that for convex bodies in $\bR^n$, of given
diameter, the maximal volume is attained for a ball. Its
reverse asks 
whether any convex body $K$ has an 
affine image $K'$ such that $(K'-K') / 2 \subset B^n$, and
$V(B^n) / V(K')$ is ``sufficiently small''? Observe that $(K'-K')/2$
is contained in $B^n$ if and only if its diameter is at most $2$,
i.e., also diam$\,K' = {\rm diam}\,[(K'-K')/2]$ is at most $2$.


{\textbf{Remark 3.20.}}
We observe that if Conjecture 3.19 holds in $\bR^3$, then
also Conjecture 3.14 holds for the case of general convex
bodies $K$ in $\bR^3$. (We note that this way was applied in
$\bR^2$, in the proof of Theorem 2 of \cite{Makai-1978}.)
In fact, consider the densest lattice packing of unit balls in
$\bR^3$. The corresponding point lattice is an inhomogeneous
lattice generated by the vertices of a regular tetrahedron. By
Theorem 1 of \cite{Makai-1978} we obtain the thinnest
non-$2$-impassable lattice of unit balls in $\bR^3$ as the polar
lattice of the densest lattice packing of $(1/4) B^3$ in
$\bR^3$.
\vskip.2cm


Recall now the considerations in the paragraph before Conjecture
3.19. Its conclusion can be formulated also in the following
way: $\min\{d_{3,2} (K) \mid K \subset \bR^3$ is a convex body$\} \ge
d_{3,2} / \min\{c > 0 \mid $ for each convex body $K \subset \bR^3$
there exists an ellipsoid $E$ such that $(K-K)/2 \subset E$ and
$V(E)/V(K) \le c\}$. The only question is whether we obtain in this
way a sharp estimate.

By Theorem 4 of \cite{Makai-1978} we have
$d_{3,2} = V (B^3) V\left( (1/4)B^3 \right) /\delta_L (B^3) = \pi/(6
\sqrt{2})$, and if the minimum of  the above numbers $c$ is as in
Conjecture 3.19, i.e., $\pi \sqrt{2}$, then we would have
$\min\{d_{3,2} (K) \mid K \subset \bR^3$ is a convex body$\} \ge
(\pi/(6\sqrt{2}))/(\pi \sqrt{2}) = 1/12$, as stated in
Conjecture 3.14 for $n = 3$.

All these point to that for finding $\min\{d_{3,2}(K) \mid K \subset
\bR^3$ is a convex body$\}$, the way through Conjecture
3.19 would be a more realistic plan than going through
Conjecture 3.17.


{\textbf{Remark 3.21.}}
Now we compare the values in Conjecture 3.14 and those following from
Theorem 3.12, first the inequalities in both chains of inequalities,
for $n=3,\dots,8,24,\,k=n-1$ (when  $d_{n,n-1}$ is known, cf., e.g.,
\cite{Rogers},
p. 3, and \cite{Cohn-Kumar1}, \cite{Cohn-Kumar2}).

For the general case the values in Conjecture 3.14 are, in the above
order,
$0,08335$..., $0,01302$..., $0,001562$..., $0,0001519$..., $0,00001240$...,
$0,0000008717$...,\\
$2,402 \cdot 10^{-30}$,
while for the centrally symmetric case they are, in the above order,
$0,1667$..., $0,04167$..., $0,008333$..., $0,001389$..., $0,0001984$...,
$0,00002480$..., $1,612 \cdot 10^{-24}$
(for $n = 24$ cf. also \cite{Gruber-Lekkerkerker}, p. 522, and \cite{Leech}.

The first inequalities in both chains of inequalities in Theorem
3.12 give, for the
general case, in the above order
$0,04538$..., $0,004548$..., $0,0003558$..., $0,00001974$...,
$0,0000008751$..., $0,00000002909$..., $4,673 \cdot 10^{-38}$,
while for the centrally symmetric case they are, in the above order,
$0,1179$..., $0,02083$..., $0,002947$..., $0,0003007$..., $0,00002482$...,
$0,000001550$..., $9,607 \cdot 10^{-32}$.\\
These mean that our estimates for $3 \le n \le 8$ can be considered as 
relatively
good.
\vskip.2cm


Now we will sharpen (\ref{equa2.11}) in \S~2. Similarly as there is a
counterpart to the question for the densest lattice packing of a
convex body $K \subset \bR^n$ (that is, of determining $\delta_L
(K)$), one can ask for an inequality from the other side:
find $\min \,\{\delta_L (K) \mid K \subset \bR^n$ is a (centrally
symmetric) convex body$\}$. This is the so-called Minkowski-Hlawka
theorem, with its variants. In our case we are interested in
$\min\,\{d_{n,k} (K)\mid K \subset \bR^n$ is a (centrally symmetric)
convex body$\}$. But here we have also a counterpart: find $\max
\{\delta_L (K) \mid K \subset \bR^n$ is a (centrally symmetric)
convex body$\}$. This is a minimax problem. Like for the
Minkowski-Hlawka theorem, for general $n$ there does not seem to
be a simple answer, what the extremal bodies would be, like for
general $n$ and $k$. For $k=0$ our question reduces to finding
$\max\,\{\delta_L (K) \mid K \subset \bR^n$ is a (centrally symmetric)
convex body$\}$. For $k=0$ and $n=2$ the solution is known,
both for the general and for the centrally symmetric case: $K$ is a
triangle, or an ellipse, respectively (cf. \cite{Gruber-Lekkerkerker},
p. 249, Theorem of I. Fáry, and p. 247, Theorem of L. Fejes Tóth,
R. P. Bambah, and C. A. Rogers).

The next interesting case is $n=2, k=1$. We improve 
(\ref{equa2.11}) in \S~2, in such a way, that the difference between the
proved and the conjectured upper bound will be reduced by a
factor about $2$.


{\textbf{Theorem 3.22.}}
\emph{We have $\max\{d_{2,1} (K) \mid K \subset \bR^2$ is a convex body$\}
= \max\{d_{2,1} (K) \mid K \subset \bR^2$ is a centrally symmetric
convex body$\} \le 0,6910 \dots$.}
\vskip.2cm


We note that the inequality of \cite{Tammela} has at its right hand side a 
certain,
explicitly given algebraic number, which is however defined in a quite
complicated way.
Therefore we do not reproduce it here. However, since \cite{Tammela}
is in Russian, we give
hints to this definition. The last displayed formula in \cite{Tammela} is the
inequality cited by us.
The quantity $m$ in it is defined by \cite{Tammela}, (18), (22), as a certain
rational function of $u_0$, where
$u_0$ is the unique root of the polynomial equation (23) in \cite{Tammela}, in
the interval given in
(23) from \cite{Tammela} (the total length of these formulas is twelve lines).


{\textbf{Conjecture 3.23.}}
Among the convex bodies $K \subset \bR^2$, the quantity
$d_{2,1} (K)$ attains its maximum for
a circle, i.e., is equal to
$\sqrt{3} \pi/8 = 0,6802...$~.
\vskip.2cm


Below (Remark 3.26) we show how R. P. Bambah's and A. C. Woods' Theorem
(see \cite{Bambah-Woods}) implies Hortob\'{a}gyi's Theorem
(cf. \cite{Hortobagyi}) (however, without the statement that there exist
even three linearly independent directions with the stated
property; i.e., his statement for one direction only). At the same
time we sharpen the Theorem of Horv\'{a}th (see \cite{Horvath-1970},
Satz 1); also in this case without the statement that there exist even four
linearly independent
directions with the stated property,
i.e., his statement for one direction only). We have to remark that this
theorem was claimed in
\cite{Horvath-1970} to
be sharp. In a paper in Hungarian
(\cite{Horvath-Ryskov-1975}, pp. 91 and 94) he has withdrawn this claim.
Below we will see that this claim does not hold; cf. our Proposition 3.27
below.
Rather than $B^n$, we will consider any convex body $K \subset
\bR^n$.


{\textbf{Theorem 3.25.}}
\emph{Let $1 \le k \le n-1$ be integers and $K \subset \bR^n$ a convex
body. Suppose $\delta_L (K) < d_{n,k} (K)$. Then for any lattice
packing $\{K+x \mid x+L \}$ of translates of $K$ there exists an
affine $k$-plane $A_k \subset \bR^n$ such that {\rm int}$(A_k
+[(d_{n,k} (K) / \delta_L (K))^{1/n}-1] (-K))$ is disjoint to
$\bigcup\{K+x\mid x \in L\}$.}
\vskip.2cm


{\textbf{Remark 3.26.}}
\textbf{A.}
Bambah-Woods' Theorem (see \cite{Bambah-Woods}) asserts $d_{3,1} = 9
\pi/32$, with equality, e.g., for the lattice generated by $(4/3)
(0,1,1),(4/3) (1,0,1),$
\newline
$(4/3) (1,1,0)$. We have $\delta_L(B^3) =
\pi/(3\sqrt{2})$ (see \cite{Rogers}, p. 3). That is smaller than
$d_{3,1} = 9 \pi/32$. Then Theorem 3.25 implies that for any
lattice packing of unit balls in $\bR^3$ there is an open,
both-way infinite cylinder with base a circle of radius
$(d_{3,1}/\delta_L(B^3))^{1/3}-1 = [(9 \pi/32) / (\pi/3
\sqrt{2})]^{1/3} - 1 = 3 \sqrt{2}/4-1$, disjoint to our lattice
packing of unit balls. That equals the value of the radius in
\cite{Hortobagyi}, Satz, cf. our (2.15) in \S~2, hence we obtained a new proof
of the inequality in \cite{Hortobagyi}, Satz. (We remind once
more that we obtained this way just one cylinder of this radius,
while \cite{Hortobagyi}, Satz, obtained three ones with linearly
independent axis directions.)

\textbf{B.} Now we show that Bambah-Woods' Theorem implies that
the value of the radius $3 \sqrt{2}/4-1$ is sharp in Hortob\'agyi's
Theorem, even when stated for one direction only --- namely for the densest
lattice packing of unit balls in $\bR^3$ --- since that was
claimed but not proved there. Suppose that there exists an open
both-way infinite cylinder of radius $r > 3 \sqrt{2}/4-1$,
disjoint to $\bigcup \{B^3 + x \mid x \in L\}$ ---
where $L$ is the point lattice
corresponding to the densest lattice packing of unit balls. Then the axis of
this
cylinder is disjoint to
$\bigcup\{(r+1){\rm int}\,B^3 + x \mid x \in L\}$. By $r+1> 3
\sqrt{2}/4$, this last lattice of open balls has a density
$\delta_L( B^3) (r+1)^3 = [\pi/(3 \sqrt{2})] (r+1)^3 >
[\pi/(3\sqrt{2})] (3 \sqrt{2}/4)^3 = 9 \pi/32 = d_{3,1}$. Now,
replacing $(r+1){\rm int}\,B^3$ by $(3 \sqrt{2}/4) B^3$, the
lattice $\{(3\sqrt{2}/4)B^3 + x\mid x \in L\}$ has a density equal to
$d_{3,1}$, and its complement contains a line. However, by
\cite{Bambah-Woods}, taking a lattice packing of closed balls of
some radius, the lattice being that of the densest packing, and
the density being $d_{3,1}$, the complement cannot contain a line.
This is a contradiction, so $r > 3 \sqrt{2}/4-1$ cannot happen,
which means that Hortob\'{a}gyi's Theorem is sharp.
\vskip.2cm



In a way very similar to that from \textbf{A}, our Theorem
3.25 implies the following proposition, that
sharpens the Theorem Horvath (\cite{Horvath-1970}, 
Satz 1), and \cite{Horvath-Ryskov-1975}, cf. our (\ref{equa2.16}) in \S~2,
although with a not explicit constant.
We remind once more that we obtain this way just
one cylinder, although with larger radius of base, 
while \cite{Horvath-1970}, Satz 1, obtained four ones with
linearly independent axis directions.

{\textbf{Proposition 3.27.}}
\emph{For each lattice packing of closed unit balls in} $\bR^4$ \emph{there
exists a both-way infinite open cylinder, disjoint to the union
of the ball lattice, with base that is a $3$-ball of radius at least
some number $ c > \sqrt{5}/2-1 \ge 0,1180\dots $\,.}


\section{Proofs}


\textbf{Proof of Proposition 3.1:} Since $d_{n,k}(K)$ is affine
invariant, we suppose that $K$ is a regular simplex $S$, say, with edge
length $\sqrt{2}$. We embed $S$ in $\bR^{n+1}$ as
conv$\{e_1,\dots,e_n,e_{n+1}\},$ where $e_i$
are the usual basis vectors. The
projection of $\{(x_1,\dots, x_n, x_{n+1}) \in \bR^{n+1} \mid \sum
^{n+1}_{i=1} x_i = 1\}$ to the $x_1,\dots,x_n$-coordinate plane
is bijective, and the vertices of $S$ project to $e_1,\dots, e_n, 0$,
respectively. The inhomogeneous lattice in $\bR^n$ generated  by
these projections is ${\Bbb Z}^n \subset \bR^n$, so the inhomogeneous lattice
generated by $e_1,\dots, e_n,e_{n+1}$ in $\bR^{n+1}$ is its inverse
image by this projection, i.e., $\{(x_1,\dots,x_n,x_{n+1}) \in
\bR^{n+1}\mid \sum ^{n+1}_{i=1} x_i = 1, \, (x_1,\dots, x_n) \in
\bZ^n\} = \{(j_1,\dots,j_n,j_{n+1}) \in \bZ^{n+1} \mid \sum
^{n+1}_{i=1} j_i = 1\}$.

Let us consider the lattice translate of the last considered
inhomogeneous lattice, i.e., $L: = \{(j_1,\dots,j_n,j_{n+1}) \in {\Bbb Z}^{n+1}
\mid \sum ^{n+1}_{i=1} j_i = 0\}$. We will determine the
D-V cell of $0$ with respect to $L$, taken in the linear hull lin$\,L$
of $L$.

Observe that the minimal length of a non-$0$
vector of $L$ is $\sqrt{2}$, and is attained exactly for $e_j - e_l,
\, 1 \le j, l\le n+1$. (Namely, any non-$0$ vector of $L$ has at
least two non-$0$ coordinates, and their absolute values are at
least 1. If its length is $\sqrt{2}$, the above two non-$0$
coordinates have absolute value 1, and all other coordinates are
$0$.)

The D-V cell of $0$ with respect to $L$, considered in lin$\,L$, is
contained in $\{x = (x_1,\dots,x_n,x_{n+1}) \in {\rm lin}\,L\mid
\|x \|^2 \le \|
x+e_j-e_l\|^2\,\, (1 \le j \not= l \le n+1)\}$. We have
\[
\begin{array}{l}
\|x\|^2 \le
\|x+e_j-e_l\|^2
\Leftrightarrow \, \langle x,x \rangle \le \langle x+e_j - e_l,
x+e_j-e_l\rangle \Leftrightarrow \\
0 \le 2 \langle x, e_j-e_l
\rangle + \langle e_j - e_l, e_j-e_l\rangle
\Leftrightarrow \, 0 \le \langle x,e_j-e_l\rangle + 1
\Leftrightarrow \\
x_l - x_j \le 1 \,.
\end{array}
\]
Hence these inequalities hold for
all $1 \le j, l \le n+1$ if and only if 
\newline
$\max _{1 \le j \le n+1}$
\newline 
$x_j - \min _{1 \le j \le n+1} x_j \le 1$. So, any $x$
in the D-V cell of $0$ with respect to $L$, considered in lin$\,L$,
satisfies the last inequality. Thus, any $x$ in the relative
interior, with respect to lin$\,L$, of this D-V cell satisfies $\max
_{1 \le j \le n+1} x_j $
\newline
$- \min _{1 \le j \le n+1} x_j
< 1$.

This also implies that any $x$ in the relative interior, with
respect to lin$\,L$, of  the D-V cell of $(i_1, \dots, i_n,i_{n+1}) \in L$
with respect to $\,L$ satisfies $\max _{1 \le j \le n+1}
(x_j$
\newline
$-i_j)- \min _{1 \le j \le n+1} (x_j - i_j) < 1$. Let $C
(i_1,\dots,i_n,i_{n+1})$, or $D(i_1,\dots,i_n,i_{n+1})$ denote the relative
interior, with respect to lin$\,L$, of the D-V cell of
$(i_1,\dots,i_n,$
\newline
$i_{n+1}) \in L$, considered in lin$\,L$, or
$\{(x_1,\dots,x_n,x_{n+1}) \in {\rm lin}\,L \mid
\max _{1 \le j \le n+1}
(x_j$
\newline
$-i_j) - \min _{1 \le j \le n+1} (x_j-i_j) < 1\}$,
respectively. Since we know that $C(0,\dots,$
\newline
$0,0) \subset
D(0,\dots,0,0)$, we have for each $(i_1,\dots,i_n,i_{n+1}) \in L$
that $C(i,\dots,i_n,$
\newline
$i_{n+1}) \subset D(i_1,\dots,i_n,i_{n+1})$.

We assert that the sets $D(i_1,\dots,i_n,i_{n+1})$, for
$(i_1,\dots,i_n,i_{n+1}) \in L$, are disjoint. Of course, it is
sufficient to show $D(0,\dots,0,0) \cap D(i_1,\dots,i_n,i_{n+1}) =
\emptyset$, for $(0,\dots,0,0) \not= \, (i_1,\dots,i_n,i_{n+1})$.
We have $\min _{1 \le j \le n+1} i_j \le -1$ or $\max
_{1 \le j \le n+1} i_j \ge 1$. Since $\sum ^{n+1}_{j=1}
i_j = 0$, any of these inequalities implies the other one (e.g.,
if $\max _{1 \le j \le n+1} i_j \ge 1$ and $\min _{1 \le j
\le n+1} i_j \ge 0$, then $\sum ^{n+1}_{j=1} i_j \ge 1$; the
other case is similar). So we have both $\min _{1 \le j \le
n+1} i_j \le -1$ and $\max _{1 \le j \le n+1} i_j \ge 1$.
Now suppose $(x_1, \dots, x_n, x_{n+1}) \in D(0,\dots,0,0) \, \cap \,
D(i_1,$
\newline
$ \dots,i_n,i_{n+1})$, and, e.g., $i_1 \le -1, i_2 \ge 1$.
Then we have $|x_1-x_2| < 1$ and $1 >$
\newline
$|(x_1-i_1)-(x_2-i_2)|=|(x_1-x_2)+(i_2-i_1)|\ge|i_2-i_1|-|x_1-x_2|>2$
\newline
$-1=1$,
a contradiction.

Thus we have $C(i_1,\dots,i_n,i_{n+1}) \subset D(i_1,\dots,i_n,i_{n+1})$
for each $(i_1,\dots,i_n,i_{n+1})$
\newline
$ \in L$, and the open
sets $D(i_1,\dots,i_n,i_{n+1})$ are disjoint for different
$(i_1,\dots,i_n$,
\newline
$i_{n+1}) \in L$. Then we
have by the above inclusion $V(C(i_1,\dots,i_n,i_{n+1})) \le
V(D(i_1,\dots$,
\newline
$i_n,i_{n+1}))$, and by the above disjointness
$V(D(i_1,\dots,i_n,i_{n+1}) \le V(C(i_1,\dots$,
\newline
$i_n,i_{n+1}))$. We
have that both $C(i_1,\dots,i_n,i_{n+1})$ and
$D(i_1,\dots,i_n,i_{n+1})$ are relative interiors of convex bodies
(in fact, polytopes) in lin$\,L$, namely of cl$\,C(i_1$,
\newline
$\dots,i_n,i_{n+1})$
and cl$\,D(i_1,\dots,i_n,i_{n+1})$. From above, the first convex
body is contained in the second one, while by the inequality
$V({\rm cl}\,C(i_1,\dots,i_n,i_{n+1})) =$
\newline
$V(C(i_1,\dots,i_n,i_{n+1})) \le V(D(i_1,\dots,i_n,i_{n+1})) = V
({\rm cl}\,D(i_1,\dots,i_n,i_{n+1}))$ the inclusion cannot be
proper. Thus, the D-V cell of $(i_1,\dots,i_n,i_{n+1})$ in $L$,
with respect to lin$\,L$, equals $\{(x_1,\dots,x_n,x_{n+1}) \in L \mid
\max _{1 \le j \le n+1} (x_j - i_j) - \min _{1 \le j \le n+1}
(x_j - i_j) \le 1\}$. In particular, the D-V cell of
$(0,\dots,0,0)$ in $L$, with respect to lin\,$L$, equals
\[
\left\{\left(x,\dots,x_n,x_{n+1}\right) \in {\rm lin}\,L \mid \max
_{1 \le
j \le n-1} x_j - \min _{1 \le j \le n+1} x_j \le
1\right\}\,.
\]
Recall the definition of $D(i_1,\dots,i_n,i_{n+1})$, in particular
that of $D(0,\dots,0,0)$ and the considerations before it. These
show that the D-V cell of $(0,\dots,0,0)$ in $L$ with respect to
lin$\,L$ equals $\{(x_1,\dots,x_n,x_{n+1}) \in {\rm lin}\, L \mid
\|x\|^2 \le
\|x+e_j-e_l \|^2 \,\, (1 \le j \not= l \le n+1)\}$. This set, $P$,
say, is a convex polytope, with at most $n(n+1)/2$ facets.
Observe, however, that $P$ is invariant under the permutations of
the coordinates $x_1,\dots,x_n,x_{n+1}$. This means that if, e.g., the
inequality $\|x\|^2 \le \|x+e_1-e_2\|^2$ would not contribute a
facet to $P$, then neither  of the inequalities $\|x\|^2 \le
\|x+e_j-e_l\|^2$ would contribute one, so $P$ would have no facets
at all, which is impossible. This shows that each inequality $\|x\|^2
\le \|x+e_j-e_l\|^2$, used for the definition of our set,
contributes a facet to $P$.

Thus the D-V cell of $(0,\dots,0,0)$ with respect to $L$, in lin$\,L$,
is a convex polytope in lin\,$L$ with $n(n+1)/2$ facets, whose affine
hulls are those affine hyperplanes in lin$\,L$, which are the perpendicular
bisectors of the segments $[0,e_j-e_l] \,\, (1 \le j \not= l \le
n+1)$. All these hyperplanes contain the  respective points
$(e_j-e_l)/2$, which lie on the boundary of the ball in lin$\,L$, with
centre $0$ and radius $1/\sqrt{2}$, and are actually tangent
hyperplanes of this ball at the respective points $(e_j-e_l)/2$.
Thus this polytope $P$, say, that evidently is a lattice space filler, is
the polar, with respect to the unit ball in lin$\,L$ with centre $0$,
of the polytope $Q:= {\mbox{conv\,}}\{e_j-e_l\mid 1 \le j \not= l \, \le n+1\}$,
that is inscribed to the ball in lin$\,L$, with centre $0$ and radius
$\sqrt{2}$.

Now observe that $Q$ is the difference body $S-S$ of the regular
simplex $S$, of edge length $\sqrt{2}$ and with
vertices $e_j \in \bR^{d+1}$, introduced at the beginning of the
proof (or we can say also with vertices $e_j-(e_1+ \dots +e_n + e_{n+1})/(n+1)
\in {\rm lin}\,L)$. We
have $V(S) = \sqrt{n+1} /n!$\,. However, we can calculate the volume
of $P$ as well. $P$ is the D-V cell of $0$ with respect to $L$,
in lin$\,L$. The orthogonal projection of $L$ to the 
$x_1\dots x_n$-coordinate plane is $\bZ^n$, that has number density $1$ in the
$x_1\dots x_n$-coordinate plane. A unit normal vector of lin$\,L$ in
$\bR^{n+1}$ is $(1/\sqrt{n+1}, \dots,
1/\sqrt{n+1}, 1/\sqrt{n+1})$. So by the above
projection the volume of the image of a set is $\langle
(1/\sqrt{n+1}, \dots, 1/\sqrt{n+1},
1\sqrt{n+1}), (0,\dots,0,1)\rangle $
\newline
$= 1/\sqrt{n+1}$ times
the volume of the original set. This means that the number density
of $L$ in lin$\,L$ is $1/\sqrt{n+1}$ times the number density of $\bZ^n$
in the $x_1 \dots x_n$-coordinate plane, i.e., it is
$1/\sqrt{n+1}$. In other words, $V(P)$, that equals the volume of a
basic parallelotope in $L$, equals $\sqrt{n+1}$.

By \cite{equa2.5},
\cite{Lovasz-Kannan}, Lemma (2.3), and \cite{Makai-1978}, Theorem 1. 
we obtain
\[
d_{n,n-1} (S) = \frac{V(((S-S)/2)^*) V (S)}{4^n \delta_L ((((S-S)/2))^*)}\,.
\]
From above we have $\delta_L (((S-S)/2)^*) = 1$. Moreover,
$V((S-S)^*)=V(Q^*)=V(P) = \sqrt{n+1}$, and $V(S)=\sqrt{n+1}/n!$.
These readily give the first equality. The second equality follows
from Stirling's formula. \hfill \rule{2mm}{2mm}


\textbf{Proof of Theorem 3.5:} Let $L \subset \bR^n$ be a lattice,
and $\lambda_1 (B^n, L) \le \dots \le \lambda_n (B^n, L)$ its
successive minima with respect to the convex body $B^n$. From
Minkowski's theorem (cf., e.g., \cite{Gruber-Handbook}, p. 750) we have
\[
\lambda_1 (B^n,L) \cdots \lambda_n (B^n,L) \le \frac{2^n
D(L)}{\kappa_n/\delta_L(B^n)} = \frac{2^n \delta_L(B^n) D(L)}{\kappa_n}\,.
\]
(Here $\kappa_n/\delta_L (B^n)$ is the smallest possible value of the
absolute value of the determinant of a packing lattice of $B^n$.)

We have
\[
\begin{array}{l}
\lambda_1 (B^n,L) \cdots \lambda_k(B^n,L) \le (\lambda_1 (B^n,L)
\cdots \lambda_n (B^n,L))^{k/n} \le\\
2^k \left(\delta_L (B^n) D(L) / \kappa_n \right)^{k/n}= 2^k
\left(\delta_L(B^n) / \kappa_n \right)^{k/n} D(L)^{k/n}\,.
\end{array}
\]
Also by \cite{Gruber-Handbook}, p. 750, we have $k$ linearly
independent vectors from $L$, with lengths $\lambda_1
(B^n,L),\dots,\lambda_k (B^n,L)$ (actually the same holds with $n$
rather than $k$). These span a $k$-dimensional sublattice, with
absolute value of determinant in
$(0,\lambda_1(B^n,L) \dots \lambda_k(B^n,L)]\subset
(0,2^k(\delta_L(B^n) / \kappa_n)^{k/n} D(L)^{k/n}]$, hence
$c_{n,k} \le 2^k (\delta_L(B^n) / \kappa_n)^{k/n}$.
The last inequality of the theorem follows from the estimate of
$\delta_L (B^n)$ in \S~1 and the formula for $\kappa_n$ (and Stirling's
formula).

For $k=1$ we have equality in the first inequality of the theorem
by the equation before the theorem, containing $c_{n,1}$.

It remained to show that for $2 \le k \le n-1$ the first
inequality of Theorem 2 is always strict. In fact, the first
inequality in the second chain of inequalities in the proof is not
sharp, unless $\lambda_1 (B^n,L) = \dots = \lambda_n (B^n,L)$.
Then the first inequality in the first chain of inequalities
reduces to the inequality
\[
\lambda_1 (B^n,L)^n \le \frac{D(L)}{2^{-n} \kappa_n/\delta_L (B^n)}\,,
\]
which is sharp only for a lattice similar to the lattice of a
densest lattice packing of $B^n$. Now consider the $n$ linearly
independent vectors from $L$, chosen with lengths $\lambda_1
(B^n,L),\dots,\lambda_n(B^n,L)$. If all of the first $k$ vectors chosen with
lengths $\lambda_1(B^n,L),\dots,\lambda_k(B^n,L)$ are not mutually
orthogonal, then the parallelepiped spanned by them has a smaller
volume than $\lambda_1(B^n,L) \cdots \lambda_k(B^n,L) $
\newline
$=\lambda_1(B^n,L)^k$. So we will have a strict inequality. Since
$\lambda_1(B^n,L)=\dots$
\newline
$=\lambda_n(B^n,L)$, we could have chosen
any $k$ of them, in any order. If not all the $n$ vectors are
mutually orthogonal, we could have chosen the first and second
ones not orthogonal, hence we have strict inequality. The same
works if among all minimal vectors of $L$ (of length
$\lambda_1(B^n,L)$) there are some not orthogonal ones. So the
only case that remains is that the lattice $L$ has only $n$
linearly independent minimum vectors, which are pairwise
orthogonal. Thus the total number of minimum vectors of $L$ is
$2n$ ($\pm 1$ times the $n$ linearly independent minimum vectors),
while in any densest lattice packing of balls the number of
minimum vectors is at least $n(n+1)$ (\cite{Gruber-Lekkerkerker},
p. 301). Hence the lattice $L$ is not similar to the lattice of
any densest lattice packing of $B^n$, a contradiction. \hfill
\rule{2mm}{2mm}


\textbf{Proof of Proposition 3.6:} Let $1 \le l \le n-k-1$ be an
integer, and let $\{e_1, \dots,e_n\}$ be a base of $\bR^n$, such that
lin$\{e_{n-k+1}, \dots, e_n\} = X$. Let the lattice $L$ have as
base $\{e_1,\dots,e_{n-k}, e_{n-k+1} + \lambda_1 e_1 + \dots +
\lambda_l e_l, e_{n-k+2}, \dots,e_n\}$, where $1, \lambda_1,\dots,
\lambda_l \in \bR$ are linearly independent over the rationals,
and let us consider the body lattice $\{x + rK\mid x \in L\}$, where
$r \in (0,\infty)$. The union of this body lattice intersects each
translate of $X$ if and only if, for the projection $\pi$  of
$\bR^n$ to ${\rm lin}\{e_1,\dots,e_{n-k}\}$ along $X$, we have that
the projection of the body lattice $\{x+rK\mid x \in L\}$ covers ${\rm
lin}\{e_1,\dots,e_{n-k}\}$.

This projection equals $P(r):= \pi(\bigcup\{i_1 e_1+ \dots +i_{n-k}
e_{n-k} + i_{n-k+1} (e_{n-k+1}+ \lambda_1 e_1 + \dots + \lambda_l
e_l) + i_{n-k+2} e_{n-k+2} + \dots + i_n e_n + rK \mid
i_1,\dots, i_n \in \bZ\} = \bigcup\{i_1e_1+\dots + i_{n-k} e_{n-k} +
i_{n-k+1} (\lambda_1 e_1 + \dots + \lambda_l e_l) + r \pi K \mid
i_1,\dots, i_{n-k+1} \in \bZ\}$. By the Theorem of Kronecker
(\cite{Bourbaki}, pp. 68-69) we have that the countably infinite
set $S:= \{i_1 e_1+ \dots+ i_l e_l + i_{n-k+1} (\lambda_1 e_1 + \dots
+ \lambda_l e_l) \mid i_1, \dots, i_l, i_{n-k+1} \in \bZ\}$ is dense in
${\rm lin} \{e_1, \dots, e_l\}$. Then $P(r) = \bigcup\{S+r \pi K +
i_{l+1} e_{l+1} + \dots + i_{n-k} e_{n-k} \mid i_{l+1}, \dots,
i_{n-k} \in {\bZ}\}$ (here $S + r \pi K$ is the Minkowski sum).
We have ${\rm lin} \{e_1, \dots, e_l\} + r \cdot {\rm rel\,int} (\pi K)
\subset S+r \pi K \subset {\rm lin}\{e_1, \dots, e_l\} + r \pi K$
(rel\,int and later rel\,bd meant with respect to the linear hull of the set
in question). Here the first or third set is a relatively open,
or closed, convex cylinder, with axis ${\rm lin}\{e_1, \dots,
e_l\}$ and base a relatively open, or closed, bounded convex set in
${\rm lin} \{e_{l+1}, \dots, e_{n-k}\}$ (the relatively
open/closed cylinder being the relative interior/closure of  the
other one). Denoting by $\varrho$ 
the projection of lin$\{e_1,\dots,e_{n-k}\}$ to
lin$\{e_{l+1}, \dots, e_{n-k}\}$, along ${\rm
lin}\{e_1,\dots,e_l\}$, we have that these bases are rel\,int$\varrho
(\pi K)$, or $\varrho (\pi K)$, respectively. Observe that $\pi K
\subset {\rm lin} \{e_1,\dots,e_{n-k}\}$ is a strictly convex
body.

The intersection of $S+r \pi K$ with the common relative
boundary of the above two cylinders is a small subset of this common
boundary. Namely, it is the union of countably many translates of
the set $[{\rm rel\,bd} (\pi K)] \cap \varrho^{-1} [{\rm rel\,bd} 
\varrho(\pi K)]$ (which is called the
\emph{shadow boundary, taken in} ${\rm lin}\{e_1, \dots,
e_{n-k}\}$, \emph{of the convex body} $\pi K \subset {\rm lin} \{e_1,
\dots, e_{n-k}\}$ \emph{with respect to illumination  from the
direction of} ${\rm lin}\{e_1, \dots, e_l\}$). The restriction of
$\varrho$ to this set is injective, since $\pi K$ is strictly convex;
so this set is topologically an $S^{d-2}$. Any of these countably
many translates is both nowhere dense in the common relative boundary of these
cylinders, and has $(n-k-1)$-Hausdorff measure $0$. There is an
$r_0 > 0$ such that $\bigcup\{r \varrho (\pi K) + i_{l+1} e_{l+1}+
\dots + i_{n-k} e_{n-k}\mid i_{l+1}, \dots, i_{n-k} \in \bZ\} = {\rm
lin} \{e_{l+1}, \dots, e_{n-k}\}$ (i.e., $\bigcup \{{\rm
lin}\{e_1,\dots,e_l\} + r \pi K + i_{l+1} e_{l+1} + \dots +
i_{n-k} e_{n-k} \mid i_{l+1}, \dots, i_{n-k} \in \bZ\} = {\rm lin}
\{e_1,\dots,e_{n-k}\})$ holds if and only if $r \ge r_0$. Then for
$r = r_0$ we have $P(r) \not= {\rm lin}\{e_1,\dots,e_{n-k}\}$.
However, for $r > r_0$ we have $P(r) \supset \bigcup \{{\rm
lin}\{e_1,\dots,e_l\}+r \cdot {\rm rel\,int} (\pi K) + i_{l+1} e_{l+1} +
\dots+ i_{n-k} e_{n-k} \mid i_{l+1}, \dots, i_{n-k} \in \bZ\} = {\rm
lin} \{e_1, \dots,e_{n-k}\}$, hence $P(r) = {\rm
lin}\{e_1,\dots,e_{n-k}\}$.

Turning back to the body lattice $\{x+rK\mid x \in L\}$, we see that its
union intersects each translate of $X$ if and only if $r > r_0$.
Dividing by $r$, we have that $\bigcup\{y+K\mid y \in L/r\}$ intersects
each translate of $K$ if and only if $r > r_0$. This is equivalent
to the claim of the proposition. \hfill \rule{2mm}{2mm}


\textbf{Proof of Theorem 3.7:}
Let $L \subset {\bR}^n$ be a lattice such that the lattice of unit balls
$\{ B^n+x \mid x \in L \} $ is $k$-impassable. We will estimate the
density of this ball lattice from below.

Let $L_k \subset L$ be a $k$-dimensional sublattice with $| {\rm det}
L_k | =D_k(L)$ and, hence, with $L \,\cap\, {\rm lin}\, L_k = L_k$.
Consider the orthogonal projection $\bR^n \to ({\rm lin}
L)^\perp$, where $\dim \,({\rm lin}\, L)^\perp = n-k$. The
image of $L$ by this projection will be a lattice, $\Lambda
\subset ({\rm lin}\, L)^\perp$, say. We have
\[
\begin{array}{c}
D(\Lambda) = D(L) / |{\rm det} L_k| = D(L)/D_k(L) \ge D(L) / (c_{n,k}
D(L)^{k/n}) =\\
D(L)^{(n-k)/n} c_{n,k} \ge D(L)^{(n-k)/n} \cdot 2^{-k}
\left(\kappa_n / \delta_L (B^n) \right)^{k/n}\,,
\end{array}
\]
where the last inequality follows from Theorem 3.4.

The projection of the $n$-dimensional lattice of unit balls will
be an $(n-k)$-dimensional lattice of unit balls in $({\rm lin}\,
L)^\perp$. Its density is $\kappa_{n-k}/D(\Lambda)$. If
$\kappa_{n-k}/D(\Lambda) < \vartheta_L (B^{n-k})$, then the projected
lattice cannot be a covering lattice for $({\rm lin}\, L)^\perp$.
Hence there is an $x \in ({\rm lin}\, L)^\perp \setminus \bigcup \{B^{n-k}
+ y \mid y \in \Lambda\}$, where in the last formula we mean by
$B^{n-k}$ the unit ball of $({\rm lin}\, L)^\perp$. Then the inverse
image of $x$ by the projection will be an affine $k$-plane in
$\bR^n$, which is a translate of lin$\,L$ disjoint to $\bigcup\{B^n +
x \mid x \in L\}$, a contradiction.

Hence we have
\[
\begin{array}{c}
\vartheta_L (B^{n-k}) \le \kappa_{n-k} / D(\Lambda) = \kappa_{n-k} D_k (L) /
D(L) \le\\
\kappa_{n-k} c_{n,k} / D(L)^{(n-k)/n} \le \kappa_{n-k} \cdot 2^k
\left(\delta_L (B^n) / \kappa_n \right)^{k/n} /
D(L)^{(n-k)/n}\,.
\end{array}
\]
Then we have for the density $\kappa_n/D(L)$ of the lattice $L \subset
\bR^n$ that
\[
\begin{array}{c}
 \kappa_n/D(L) \ge 
\frac{\kappa_n \vartheta_L (B^{n-k})^{n/(n-k)}}
{\kappa_{n-k}^{n/(n-k)} c_{n,k}^{n/(n-k)}} \ge \\
\frac{\kappa_n \vartheta_L (B^{n-k})^{n/(n-k)}}{\kappa_{n-k}^{n/(n-k)}
2^{kn/(n-k)}
\left(\delta_L(B^n) / \kappa_n \right)^{k/(n-k)}}\,.
\end{array}
\]
Since the first inequality of Theorem 3.4 was strict for $k \ge
2$, also the se\-cond inequality of Theorem 3.5 is strict for $k \ge
2$.

On the other hand, we have $d_{n,k} \ge d_{n,n-1} =
k^2_n/(4^n \delta_L (B^{(n)}) = n^{-n} e^{O(n)}$ by
\S~1, (\ref{equa2.6}) in \S~2, and by Stirling's formula. \hfill
\rule{2mm}{2mm}


\textbf{Proof of Corollary 3.9:} We evaluate the second lower
bound from Theorem 3.7 for $n = 2$, using $\delta_L (B^2)
= \pi/\sqrt{12}$ and $\vartheta_L (B^1) = 1$.

The sharpness of the estimate for $n=2$ follows from
\cite{Makai-1978}. \hfill\rule{2mm}{2mm}


\textbf{Proof of Proposition 3.11:} By (2.18) from the
introduction there exists a constant $n_0$ such that for any $n
\ge n_0$ and $n/ \log _2 n \le k \le n-1$ there exists a lattice packing
$\{B^n + x\mid x \in L\}$ of unit balls that is $k$-impassable. Then
$d_{n,k}$ is at most the density of this lattice packing, which is
in turn at most $\delta_L (B^n)$. This proves the first inequality.
For the second inequality cf. in the preliminaries. \hfill \rule{2mm}{2mm}


\textbf{Proof of Theorem 3.12:} By \cite{Ball},
\cite{Barthe-1997}, \cite{Barthe-1998} (the paragraph before
Corollary 3, and Proposition 10), and
\cite{Barthe-2003} (pp.
147-148) a convex body, or a centrally symmetric convex body $K
\subset \bR^n$ can be included into an ellipsoid $E$ such that
$V(K)/V(E) \ge V(S^n)/V(B^n)=(n+1)^{(n+1)/2}/(\kappa_n n!n^{n/2})$ or
$V(K)/V(E)\ge V(B^n)/V(C^n)=2^n/(\kappa_n n!)$, where $S^n$, or $C^n$,
are the regular simplex, or regular cross-polytope, inscribed to
$B^n$, respectively. (Furthermore, the only cases of equality are
that $K$ is a simplex, or cross-polytope, respectively.)

Now let us include our convex body, or centrally symmetric convex
body $K$, into such an ellipsoid $E$,
as stated above. If for a point
lattice $L \subset \bR^d$ we have that the body lattice $\{K+x\mid x
\in L\}$ is $k$-impassable, then also the body lattice $\{E+x\mid x
\in L\}$ is $k$-impassable, with the quotient of the density of
the mentioned lattice of translates of $K$ and the density of the
mentioned lattice of translates of $E$ being $V(K) / V(E)$ (for
which quantity we have the lower estimates given above). The second body
lattice is a lattice
of ellipsoids. If $E$ were a ball, the density of this second body
lattice would be at least $d_{n,k}$, by definition. However,
$d_{n,k}(K)$ is invariant under affine transformations of
$K$. So the density of the second body lattice is at least
$d_{n,k}$, anyway.

Then the density of the first body lattice (of translates of $K$)
is the product of the density of the second body lattice (for
which the lower estimate $d_{n,k}$ holds) and the quotient of
the densities of the first body lattice and the second body
lattice (that equals $V(K)/V(E)$, for which we have the above
given lower estimates). The product of these lower estimates gives
the first inequality (both for the general and the centrally
symmetric case). The second and third inequalities follow from
Theorem 3.7.

On the other hand, we have $d_{n,k} (K) \ge d_{n,n-1} (K) >
\kappa _n^2 {2n \choose n}^{-1}=e^{O(n)}/n^n$, by (\ref{equa2.10}).
\hfill \rule{2mm}{2mm}


\textbf{Proof of Theorem 3.16:} \textbf{A.}
By \cite{Reisner}, for $K$ a zonoid
centred at $0$, or a polar of a zonoid centred at $0$, we have
that Mahler's conjecture is true, i.e., $V(K) V(K^*) \ge 4^n/n!$.
Then (\ref{equa2.5}) gives
\[
d_{n,n-1} (K) = V (K) V(K^*) / [4^n \delta_L (K^*)] \ge V(K)
V(K^*) / 4^n \ge 1/n!\,.
\]
For $K$ a cross-polytope we have that $K^*$ is a parallelotope,
hence $\delta_L (K^*) = 1$, and so the second minimum
in the theorem equals $1/n!$.

By \cite{Reisner}, for $K$ a zonoid centred at $0$ we have $V(K)V(K^*) >
4^n/n!$,
unless $K$ is a parallelotope. For $K$ a parallelotope centred at
$0$ we have that $K^*$ is a cross-polytope. Then $\delta_L (K^*) <
1$, since else a lattice of translates of a cross-polytope would tile
$\bR^n$, and so all its facets ought to be centrally symmetric (see
\cite{Gruber-Lekkerkerker}, p. 168), which is false for $n \ge 3$. Now observe
that the expression for $d_{n,n-1}(K)$ is affine invariant and
continuous in $K$, hence by \cite{Schneider}, p. 60, Note 13, the
first minimum in Theorem 3.16 exists (is attained). The same
consideration proves that the second minimum in the theorem is attained
only for a cross-polytope.

\textbf{B.}
Analogously, for the bodies in the third and fourth minima, Mahler's
conjecture is true, by \cite{Saint-Raymond}, Th\'eor\`eme 25 and 28.
Then an analogous consideration
proves that this minimum is also equal to $1/n!$.

We turn to show that, in the third minimum, the minimum is attained only for
cross-polytopes. For this aim, we have to recall a result of \cite{Meyer},
Th\'eor\`eme 1.3 and 1.4, that establishes all cases of equality in the
inequality $V(K)V(K^*) \ge 4^n/n!$, for $K \subset \bR^n$ a convex body
symmetric to all coordinate hyperplanes. All these
cases of equality are obtained in
the following way. Let the dimension $n \ge 1$ be fixed. We consider $\bR^n$,
with the standard base $\{ e_1,...,e_n \}$. We will work with coordinate
subspaces of $\bR^n$. We begin with defining $0$-symmetric convex bodies in
all $1$-dimensional coordinate subspaces of $\bR^n$, namely the segments
$[-e_i,e_i]$. We proceed by induction. Let us have a set of coordinate
subspaces $\{ X_j \}$ in $\bR^n$, with pairwise intersections $\{ 0 \}$,
and together spanning $\bR^n$, and in each $X_j$ let us have an $0$-symmetric
convex body $K_j$. If there are at least two such $X_j$-s, we pick
two of them, and replace this pair by their (direct) sum, and replace the
corresponding two $K_j$-s either by their (direct) sum, or by the convex hull
of their union. We end when we have only one $X_j$, and the corresponding
$K_j$ will be the body constructed this way. (The normed spaces
corresponding to these convex bodies are called
\emph{Hanner-Hansen-Lima spaces}, because these authors 
investigated their properties.)
Clearly, all bodies constructed this way are polytopes
symmetric w.r.t. all coordinate
hyperplanes. Furthermore, the polar of such a polytope is such a polytope
as well.

We are going to show that, unless $K$ is a cross-polytope (i.e., $K^*$ is a
parallelotope), we have $\delta _L (K^*)<1$. For this, like in \textbf{A},
it will suffice to show that $K^*$ has a facet which is not centrally 
symmetric.
We will use induction for $n$.

We will use that also $B:=K^*$
is obtained by the above construction. In the last
step of the construction, we will have two coordinate subspaces $X_1,X_2$, and
convex polytopes $B_1,B_2$ in them.

If $B$ is the direct sum of $B_1$ and
$B_2$, and is not a parallelotope, then one $B_j$ is not a
parallelotope. Then, by the induction hypothesis, $B_j$
has a facet $F_j$, say, which is not centrally symmetric. Then $B$ has a facet 
$F_j
\oplus B_{2-j}$ that is not centrally symmetric either.

Now let $B$ be the convex hull of the union of $B_1$ and $B_2$. By turning to
the polar bodies, we see that the vertices of $B^*$ are the direct sums of
any vertex of $B_1^*$ and any vertex of
$B_2^*$. Hence, the facets of $B$ are the convex
hulls of the unions of any facets of $B_1$ and $B_2$. Now recall that
$n \ge 3$. (We
remark that the case $n=2$ is anyway trivial, since the direct sum and the
convex hull of the union of $B_1$ and $B_2$ both are parallelograms.)

Let therefore $F_j$ be a facet of $B_j$, and consider the facet $F=$conv$\,(F_1
\cup F_2)$ of $B$. Observe that $2 \le n-1={\rm dim}\,F
= {\rm dim}\,F_1+{\rm dim}\,F_2+1$. Hence one of dim$\,F_j$ is positive.
Let us consider aff$\,F$. Observe that $({\rm aff}\,F_1) \cap
({\rm aff}\,F_2)$
is a subset of the intersection of the corresponding coordinate
hyperplanes, i.e., of $\{ 0 \} $. However, $F_j \ni 0$, so $({\rm aff}\,F_1)
\cap ({\rm aff}\,F_2)$ is empty. There is, up to translations,
just one hyperplane $H$ in aff$\,F$
that is parallel to both aff$\,F_j$. The boundary hyperplanes of the
supporting strip of $F$
in aff$\,F$, parallel to $H$, intersect $F$
in $F_1$ and $F_2$, respectively. If $F$ would be centrally symmetric,
the affine hulls of these intersections
whould be translates of each other. So, if we consider their translates
containing $0$, these would coincide. However, these translates have
intersection $\{ 0 \} $. So these translates would be $\{ 0 \} $, so both
would have dimension $0$. This, however, contradicts the fact that
one of dim$\,F_j$ is positive.

\textbf{C.}
Lastly, for $n \le 8$, also for the
bodies in the fourth and fifth minima Mahler's
conjecture is true, by \cite{Lopez-Reisner}. Even, they prove that the
only cases of equality are of the form described in \textbf{B}. Then the
considerations of \textbf{B} show that the only case, when in
\[
d_{n,n-1} (K) = V (K) V(K^*) / [4^n \delta_L (K^*)] \ge V(K)
V(K^*) / 4^n \ge 1/n!
\]
we have equality at both inequalities, is when $K$ is a
cross-polytope.

This shows that the fifth minimum in the theorem is greater
than $1/n!$, and the sixth minimum is $1/n!$, and is attained only for $K$ a
cross-polytope.
\hfill
\rule{2mm}{2mm}


\textbf{Proof of Theorem 3.18:} The considerations after Conjecture
3.17 give, for the minimum in
the theorem, the lower estimate
\[
\begin{array}{c}
\left[ 2^n \big/ {2n \choose n} \right] \cdot \min \{ V(K)V(K^*) \mid K
\subset \bR^n \\
\mbox{is an 0-symmetric convex body} \}\,.
\end{array}
\]
Then we apply the theorem of Reisner \cite{Reisner}
used in the proof of Theorem 3.16.

For $K$ a simplex, Proposition 3.1 gives $d_{n,n-1} (K)
= 2^n (n+1)/n!$. It remains to show that, for $K$ a simplex,
$((K-K)/2)^*$ is a zonoid. We are going to show that it is a
zonotope.

Like in the proof of Proposition 3.1, we may suppose that $K$
equals a regular simplex with edge length $\sqrt{2}$, which we
denote by $S$. Further we use the notations from the proof of
Proposition 3.1. There it was proved that the D-V
cell of $(0,\dots,0,0)$ with respect to $L$, in ${\rm lin}\,L$, is the polar
of the polytope $Q = {\rm conv} \{e_j - e_l \mid 1 \le j \not= l
\le n+1\} = S-S$, with respect to the unit ball in lin$\,L$ with
centre $0$. So the above D-V cell is $(S-S)^*$. By \cite{Vegh},
this D-V cell is the orthogonal projection of the D-V cell of
$(0,\dots,0,0)$ with respect to $\bZ^{n+1}$, in $\bR^{n+1}$, i.e.,
of the unit cube $C=[-1/2,1/2]^{n+1}$. (Considering the nearest
neighbours of $(0,\dots,0,0)$ in $\bZ^n$ gives that this D-V cell
is contained in $C$ and also has volume $1$, so they are equal.)
Here $C$ is a zonotope, hence its orthogonal projection $(S-S)^*$ is a
zonotope, too. And $((S-S)/2)^* = 2 (S-S)^*$ is a zonotope as well.
\hfill \rule{2mm}{2mm}


\textbf{Proof of Theorem 3.22:} The equality follows from the
inequality $d_{2,1} (K) \le d_{2,1} ((K-K)/2)$, cf.
\cite{Makai-1978}, Theorem 4. Now we prove the inequality for $K$
centrally symmetric.

By \cite{Makai-1978}, Theorem 5, we have, for centrally symmetric
$K$,
\[
d_{2,1} (K) = \frac{\pi^2}{16 \min \{\delta_L(K')\}}\,,
\]
where the minimum is taken for all centrally symmetric convex
bodies $K' \subset \bR^2$. This minimum satisfies
\[
\min \{\delta_L (K')\} \ge 0.8926\dots\,,
\]
cf. \cite{Tammela}. These two inequalities imply our theorem.
\hfill \rule{2mm}{2mm}


\textbf{Proof of Theorem 3.25:} The density of the body lattice
$\{K+x\mid x \in L\}$ is at most $\delta_L (K)$, that is smaller than
$d_{n,k}(K)$. Hence there exists an affine $k$-plane disjoint to
$\bigcup \{K+x\mid x\in L\}$. Even, we may inflate $K$ by a factor
$(d_{n,k} (K) / \delta_L(K))^{1/n}$, and still we find an
affine $k$-plane $A_k$ disjoint to $\bigcup\{{\rm int}\, (d_{n,k} (K) /
\delta_L(K))^{1/n} K+x\mid x\in L\}$. Inflation of $K$ by a factor
$(d_{n,k} (K) / \delta_L (K))^{1/n}$ can be written also as
$K+[(d_{n,k}(K)/\delta_L(K))^{1/n}-1]K$.
Then disjointness of $\{{\rm int} (K+[(d_{n,k} (K)/ \delta_L
(K))^{1/n}-1]K)+x\mid x \in L\}$ and $A_k$ is equivalent to
disjointness of $\{K+x\mid x \in L\}$ and int$\,(A_k+[(d_{n,k} (K) /
\delta_L(K))^{1/n}-1] (-K))$. \hfill \rule{2mm}{2mm}


\textbf{Proof of Proposition 3.27:.}
By Theorem
3.7, we have $d_{4,1} > 25 \pi^2/256$. We have $\delta_L (B^4)
= \pi^2/16$ (\cite{Rogers}, p. 3), which is smaller than $25 \pi^2/256 \,\,(\le
d_{4,1})$. Then,
by Theorem 3.25, applied to $B^4 \subset \bR^4$, there
is an open, both-way infinite cylinder with base a $3$-ball of
radius $(d_{4,1}/\delta_L(B^4))^{1/4}-1 >
[(25\pi^2/256)/(\pi^2/16)]^{1/4}-1 = \sqrt{5}/2-1$, which is disjoint to
our lattice packing of closed unit balls. 
$\hfill \rule{2mm}{2mm} $


\vspace{3cm}

E. Makai, Jr.\\
A. Rényi Institute of Mathematics\\
Hungarian Academy of Sciences\\
PF 127\\
H-1364 Budapest\\
HUNGARY\\
makai@renyi.hu\\[0.3cm]

H. Martini\\
Mathematical Faculty\\
University of Technology\\
D-09107 Chemnitz\\
GERMANY\\
martini@mathematik.tu-chemnitz.de

\end{document}